\documentclass[preprint,12pt]{elsarticle}

\usepackage{amsmath,amsthm,amssymb,multirow,mathtools}

\usepackage[all]{xy}
\usepackage[bookmarks,colorlinks]{hyperref}

\usepackage{color}
\usepackage[linesnumbered,ruled,vlined]{algorithm2e}
\usepackage{tikz}
\usepackage{tikz-cd}

\numberwithin{equation}{section}

\newtheorem{theorem}{Theorem}

\theoremstyle{definition}
\newtheorem{definition}{Definition}[section]
\newtheorem{lemma}[definition]{Lemma}%Bổ đề
\newtheorem{proposition}[definition]{Proposition}%Mệnh đề
\newtheorem{corollary}[definition]{Corollary}%Hệ quả
\newtheorem{example}[definition]{Example}%{Ví dụ}
\newtheorem{remark}[definition]{Remark}%{nhận xét}

\newcommand{\A}{\mathcal{A}}
\newcommand{\ad}{\mathrm{ad}}
\newcommand{\aff}{\mathfrak{aff}}
\newcommand{\Au}{\mathrm{Aut}}
\newcommand{\B}{\mathcal B}
\newcommand{\C}{\mathbb C}
\newcommand{\Cc}{\mathcal C}
\newcommand{\co}{\coloneqq}
\newcommand{\Der}{\mathrm{Der}}
\newcommand{\GL}{\mathrm{GL}}
\newcommand{\id}{\mathrm{id}}
\newcommand{\im}{\mathrm{im}\,}
\newcommand{\K}{\mathbb K}
\newcommand{\N}{\mathbb N}
\newcommand{\M}{\mathrm{Mat}}
\newcommand{\R}{\mathbb R}
\newcommand{\rc}{\mathcal R}

\newcommand{\s}{\mathrm{span}}
\newcommand{\Sc}{\mathcal S}
\newcommand{\Sf}{\mathfrak s}
\newcommand{\Z}{\mathcal Z}
\newcommand{\LO}{\overset{\rm Lie}{\oplus}}

\begin{document}

\begin{frontmatter}

\title{Complete and cocomplete Lie algebras with injective- and projective-type properties}

\author[1,2]{Vu A. Le}
\ead{vula@uel.edu.vn}
\address[1]{University of Economics and Law, Ho Chi Minh City, Vietnam}
\address[2]{Vietnam National University, Ho Chi Minh City, Vietnam}

\author[3]{Hoa Q. Duong}
\ead{hoa.duongquang@ufm.edu.vn}
\address[3]{University of Finance - Marketing, Ho Chi Minh City, Vietnam}

\author[4]{Tuan A. Nguyen\corref{cor1}}
\ead{tuannguyenanh@hcmue.edu.vn}
\address[4]{Ho Chi Minh City University of Education, Ho Chi Minh City, Vietnam}
\cortext[cor1]{Corresponding author: tuannguyenanh@hcmue.edu.vn}

%-------Abstract -----
\begin{abstract}
In the category of modules, injective and projective objects are characterized by the splitting of short exact sequences. Motivated by this principle, we investigate analogous phenomena in the category of finite-dimensional Lie algebras over a field of characteristic zero. Since this category is not abelian, extensions admit two distinct notions of splitting, which we call \emph{trivial} and \emph{semi-trivial}. The first main result establishes the converse of Jacobson's classical theorem: a Lie algebra $\A$ trivially splits every extension by it if and only if $\A$ is complete. This identifies completeness as the ``injective-type'' property in this category -- although, seemingly weaker than category-specific injectivity. By contrast, no nontrivial Lie algebra has the dual property: the second main result asserts that, for every nontrivial Lie algebra $\Cc$, some extension of $\Cc$ fails to split trivially, so no ``projective-type'' analogue exists. Restricting to central extensions restores a workable dual
notion: a Lie algebra $\Cc$ is called cocomplete if every central extension of it splits trivially, and we prove that this holds if and only if $H^2(\Cc, \K) = 0$; in
particular, semisimple Lie algebras are both complete and cocomplete. Each of these results admits an equivalent homomorphism-lifting reformulation, paralleling the lifting properties of injective and projective modules. For
almost abelian Lie algebras, cocompleteness reduces to an explicit spectral condition on the defining derivation. These characterizations underlie three corresponding algorithms, which yield tabulations of the complete, the cocomplete, and the almost abelian cocomplete Lie algebras of dimension at most 4. 
\end{abstract}
\end{frontmatter}

%-------Section 1: Introduction ----
\section{Introduction}\label{sec1}
A fundamental problem in many algebraic categories is to determine whether a short exact sequence splits. In the category of modules, projective and injective objects are characterized by the trivially split property of short exact sequences. Specifically, let $R$ be a unital ring and

%----A Extension of Modules-----
\begin{equation}\label{ME}
    \xymatrix
	{
 		0 \ar[r] & A \ar[r]^\iota & B \ar[r]^\pi & C \ar[r] & 0
 	}\tag{ME}
\end{equation}
is a short exact sequence of right $R$-modules. By Jacobson~\cite[Lemma~3.9.3]{Jac62}, the following three statements are equivalent.

%-----A Extension of Modules splits (trivially)---
\begin{enumerate}
    \item[(i)] The short exact sequence~\eqref{ME} splits trivially, i.e., $B \cong A \oplus C$ as a direct sum of right $R$-modules.
    \item[(ii)] The module homomorphism $\iota$ admits a left inverse, i.e, there exists a module homomorphism $j \colon B \to A$ with $j \circ \iota = \id_A$.
    \item[(iii)] The homomorphism $\pi$ admits a right inverse, i.e., there exists a module homomorphism $q \colon C \to B$ with $\pi \circ q = \id_C$.
\end{enumerate}
Furthermore, the right $R$-module $A$ is injective if and only if every short exact sequence~\eqref{ME} splits trivially. Similarly, the right $R$-module $C$ is projective if and only if every short exact sequence~\eqref{ME} splits trivially.

This naturally raises the question whether a similar phenomenon occurs in the category of Lie algebras. Consider a short exact sequence of Lie algebras
%---Lie algebra Extension------
\begin{equation}\label{LE0}
    \xymatrix
	{
 		0 \ar[r] & \A \ar[r]^\iota & \B \ar[r]^\pi & \Cc \ar[r] & 0
 	}\tag{LE}
\end{equation}
which is called a \emph{Lie algebra extension} (for short, \emph{extension}) of $\Cc$ by $\A$. A natural problem is to find structural conditions on $\A$ or $\Cc$ ensuring that extension~\eqref{LE0} splits trivially, i.e., $\B \cong \A \LO \Cc$.

However, the category of Lie algebras differs from that of modules in a fundamental way. By Levi's theorem~\cite[Part~I, Chapter~VI, Theorem~4.1]{Serre}, whenever $\Cc$ is semisimple, the epimorphism $\pi$ in extension~\eqref{LE0} always admits a Lie algebra right inverse $q \colon \Cc \to \B$ with $\pi \circ q = \id_{\Cc}$; nevertheless, this generally realises $\B$ only as a semidirect, not direct, sum of Lie algebras. We therefore distinguish two notions of splitting for an extension~\eqref{LE0}:

%----Semitriviall splitting/Trivial splitting-----
\begin{itemize}
    \item \emph{Semi-trivial splitting}: when $\B \cong \A \oplus_{\rho} \Cc$ for some action $\rho$ of $\Cc$ on $\A$;
    \item \emph{Trivial splitting}: when $\B \cong \A \LO \Cc$.
\end{itemize} 
These notions give rise to the following two basic questions:
\begin{itemize}
\item Under which conditions on $\A$ do all extensions of $\Cc$ by $\A$ split trivially?
\item Under which conditions on $\Cc$ do all extensions of $\Cc$ by $\A$ split trivially?
\end{itemize}

Regarding the first question, a motivation comes from a classical result of Jacobson~\cite[Proposition~3]{Jac62} which asserts that if $\A$ is a complete Lie algebra (i.e., $\A$ has trivial center and only inner derivations), %(see Definition~\ref{D34}), 
then every extension~\eqref{LE0} is trivially split. The first main contribution of the paper is to prove that the converse also holds: for any Lie algebra $\A$, the trivial splitting of every extension~\eqref{LE0} is equivalent to the condition that $\A$ is complete. In other words, we prove that \emph{complete Lie algebras have an ``injective-type'' property within the category of Lie algebras}, although, as Example~\ref{Ex-counter-cat-inj} shows, this property is strictly weaker than categorical injectivity. This result is formulated as Theorem~\ref{thm1} below.

%-------Theorem 1 (New)------
\begin{theorem}[Injective-type characterizations of complete Lie algebras]\label{thm1}
For a Lie algebra $\A$, the following statements are equivalent.
\begin{enumerate}[(1)]
    \item\label{Thm1-1} $\A$ is complete.
    \item\label{Thm1-2} \emph{(Lifting homomorphisms into $\A$)} Suppose that the diagram
    \begin{equation}\label{Diag1}
%        \begin{tikzcd}[column sep=large]
%            0 \arrow[r] & \A' \arrow[r, "f"] \arrow[d, "g"] & \B \\ 
%            & \A &
%        \end{tikzcd}
        \xymatrix
	       {
 		 0 \ar[r] & \ar[d]_g \A' \ar[r]^f & B \\
            & \A
 	      }
    \end{equation}
    satisfies the following assumptions:
    \begin{enumerate}[(2i)]
        \item\label{Thm1-2i} The horizontal line is exact, i.e., $f$ is a Lie algebra monomorphism, with $\im f$ is an ideal of $\B$.
        \item\label{Thm1-2ii} $\im g$ is a direct summand of $\A$, i.e., $\A = \im g \LO J$ for some ideal $J$ of $\A$.
        \item\label{Thm1-2iii} $f(\ker g)$ is an ideal of $\B$.
    \end{enumerate}
    Then there exists a Lie algebra homomorphism $h \colon \B \to \A$ with $\im h = \im g$ such that $h \circ f = g$, i.e., the diagram~\eqref{Diag1} always embeds in the following commutative diagram.  
    \[
%        \begin{tikzcd}[column sep=large]
%            0 \arrow[r] & \A' \arrow[r, "f"] \arrow[d, "g"] & \B \arrow[ld, dashed, "\exists\,h"] \\
%            & \A &
%        \end{tikzcd}
        \xymatrix
	{
 		0 \ar[r] & \ar[d]_g \A' \ar[r]^f & B \ar@{-->}[dl]^{\exists \, h}  \\
        & \A
 	}
    \]
    \item\label{Thm1-3} Every extension~\eqref{LE0} of any Lie algebra $\Cc$ by $\A$ splits trivially.
    \item\label{Thm1-4} For every Lie algebra monomorphism $\iota \colon \A \hookrightarrow \B$ whose image $\im\iota$ is an ideal of $\B$, $\im\iota$ is always a direct summand of $\B$, %as a Lie algebra; 
    i.e., $\B = \im\iota \LO I$ for some ideal $I$ of $\B$.
\end{enumerate}
\end{theorem}

Regarding the second question, a dual inquiry arises: does there exist a class of Lie algebras that play a ``projective-type'' role, mirroring the behavior of projective modules in Module Theory? Although Theorem~\ref{thm1} provides a complete characterization of the injective-type property, the projective counterpart turns out to be elusive. The second main result reveals that such a property is too restrictive to exist globally in the category of Lie algebras.

%----Theorem 2 (New)------------
\begin{theorem}[{Nonexistence of projective objects in the category of Lie algebras}]\label{thm2}
No nontrivial Lie algebra possesses a global projective-type property. More precisely, for every nontrivial Lie algebra $\Cc$, there exist a Lie algebra $\A$ and an extension~\eqref{LE0} of $\Cc$ by $\A$ that does not split trivially, i.e., $\B \ncong \A \LO \Cc$.
\end{theorem}
	
This nonexistence reveals a structural asymmetry between the injective and projective aspects of Lie algebra extensions.
It suggests that the projective-type property cannot be realized in the absolute sense, but rather must instead be sought within a more restricted framework. A natural choice is the class of \emph{central extensions}~\cite[Definition 4.1]{Sch08}% (see Definition~\ref{DCLE} in Section~\ref{sec2} or~\cite[Definition 4.1]{Sch08}) 
\begin{equation}\label{CE0}
    \xymatrix
	{
	   0 \ar[r] & V \ar[r] & \B \ar[r] & \Cc \ar[r] & 0
	}\tag{CE}
\end{equation}
in which $\im \iota = \ker \pi$ is contained in the center $\Z(\B)$ of $\B$. The motivation for this restriction is rooted in classical Lie theory: a well-known consequence of Levi's theorem is that every central extension of a semisimple Lie algebra splits trivially (Proposition~\ref{P312}), so semisimple Lie algebras already exhibit a projective-type behaviour with respect to central extensions.

However, semisimplicity is sufficient but not necessary, and the property ``every central extension splits trivially'' therefore deserves to be treated as a distinct structural condition. We accordingly introduce a new class of Lie algebras singled out by this property, which we call \emph{cocomplete} Lie algebras (see Definition~\ref{D311}) since the name reflects their duality to complete Lie algebras. The third main result in Theorem~\ref{P-lifting-cocomplete} below, formally dual to Theorem~\ref{thm1}, presents several equivalent characterizations of cocompleteness, paralleling those for complete Lie algebras and matching a familiar property in the category of modules (compare~\cite[Theorem~5.3.1]{Kas82}).

%---Theorem 3 (New)-------------
\begin{theorem}[Projective-type characterizations with respect to central extensions]\label{P-lifting-cocomplete}
For a Lie algebra $\Cc$, the following statements are equivalent.
\begin{enumerate}[(1)]
    \item\label{Thm3-0} $\Cc$ is cocomplete.
    \item\label{Thm3-1} Every central extension \eqref{CE1} of $\Cc$ by any vector space $V$ splits trivially.
    \item\label{Thm3-2} For every Lie algebra epimorphism $\pi \colon \B \twoheadrightarrow \Cc$ with $\ker\pi \subseteq \Z(\B)$, the kernel $\ker\pi$ is a direct summand of $\B$, i.e., $\B = \ker\pi \LO I$ for some ideal $I \cong \Cc$ of $\B$.
    \item\label{Thm3-3} \emph{(Lifting homomorphisms from $\Cc$)} Suppose that the diagram
    \begin{equation}\label{Diag2}
%        \begin{tikzcd}
%            & \Cc \arrow[d, "g"] \\ %\arrow[ld, dashed, "\exists\,h"'] & \\
%            \B \arrow[r, "f"'] & \Cc' \arrow[r] & 0
%        \end{tikzcd}
        \xymatrix
	{
	   & \Cc \ar[d]^g \\
       \B \ar[r]_f & \Cc' \ar[r] & 0
	}
    \end{equation}
    satisfies the following assumptions:
    \begin{enumerate}[(4i)]
        \item\label{Thm3-3i} The horizontal line is exact, i.e., $f$ is a Lie algebra epimorphism. % with $\im f$ is an ideal of $\B$.
        \item\label{Thm3-3ii} $\ker f \subseteq \Z(\B)$. 
    \end{enumerate}
    Then there exists a Lie algebra homomorphism $h \colon \Cc \to \B$ such that $f \circ h = g$, i.e., the diagram~\eqref{Diag2} always embeds in the following commutative diagram.
    \[
%        \begin{tikzcd}
%            & \Cc \arrow[d, "g"] \arrow[ld, dashed, "\exists\,h"'] & \\
%            \B \arrow[r, "f"'] & \Cc' \arrow[r] & 0
%        \end{tikzcd}
         \xymatrix
	{
	   & \ar@{-->}[dl]_{\exists \, h} \Cc \ar[d]^g \\
      \B \ar[r]_f & \Cc' \ar[r] & 0
	}
    \]
\end{enumerate}
\end{theorem}

After proving Theorems~\ref{thm1},~\ref{thm2} and ~\ref{P-lifting-cocomplete}, the remainder of the paper is devoted to the fundamental properties of cocomplete Lie algebras. By cohomological methods, we obtain an explicit characterization: a Lie algebra $\Cc$ is cocomplete if and only if $H^2(\Cc,\K)=0$. In particular, every semisimple Lie algebra is cocomplete. Thus, the class of semisimple Lie algebras lies in the intersection of the complete and the cocomplete ones.

Although cocomplete Lie algebras possess several interesting properties, their classification is nontrivial. A simple observation shows that many low-dimensional cocomplete Lie algebras are \emph{almost abelian}, which motivates our study of Lie algebras that are both almost abelian and cocomplete. Recall that an $(n+1)$-dimensional non-abelian Lie algebra $\Cc$ is almost abelian if it admits an $n$-dimensional abelian subalgebra~\cite[Definition~1]{Ave22}, and such a subalgebra is automatically an ideal~\cite[Proposition 3.1]{BC12}, so there exists an abelian ideal $\K^n \subset \Cc$ such that $[\Cc, \K^n] \subset \K^n$. Equivalently, $\Cc$ is the 1-dimensional extension of $\K^n$ by a derivation $D \in \Der(\K^n) = \mathfrak{gl}_n(\K)$, namely,
\[
	\Cc = \K^n \oplus_D \K e_0, \quad e_0 \in \Cc \setminus \K^n.
\]
Within this framework, the final main result gives a necessary and sufficient condition for an almost abelian Lie algebra to be cocomplete as follows.

%--------Theorem 4-------
\begin{theorem}[{Cocompleteness criterion for almost abelian Lie algebras}]\label{thm3}
    Let $\K$ be a field of characteristic $0$ with the algebraic closure $\bar{\K}$. An $(n+1)$-dimensional almost abelian Lie algebra $\Cc = \K^n \oplus_D \K e_0$ is cocomplete if and only if the following conditions hold.
    \begin{enumerate}[(1)]
	   \item\label{thm3-1} $D$ is an automorphism of $\K^n$, i.e., $D \in \Au(\K^n) \equiv \GL_n(\K)$. 
	   \item\label{thm3-2} No pair of eigenvalues $\lambda$, $\mu$ of $D_{\bar{\K}} \co D \otimes_\K \id_{\bar{\K}}$ satisfies $\lambda + \mu = 0$.
    \end{enumerate}
\end{theorem}

Theorem~\ref{thm3} reduces the classification of almost abelian cocomplete Lie algebras to that of invertible matrices satisfying condition~\eqref{thm3-2}, up to \emph{proportional similarity} (see Proposition~\ref{P354}). Combining with computer algebra techniques in \cite{NLV25}, this yields a classification of low-dimensional almost abelian cocomplete Lie algebras up to isomorphism.

The paper is organized as follows. Section~\ref{sec2} recalls Lie algebra extensions, in particular the central ones. Section~\ref{sec3} is devoted to complete Lie algebras in which we present their basic properties, prove Theorem~\ref{thm1} and conclude with a practical algorithm for verifying completeness. Section~\ref{sec4} is the parallel section for cocomplete Lie algebras: we prove Theorems~\ref{thm2} and~\ref{P-lifting-cocomplete}, give the cohomological characterization, provide a practical algorithm for testing cocompleteness, and finally prove Theorem~\ref{thm3} together with a classification algorithm for almost abelian cocomplete Lie algebras.

\paragraph{Notation and conventions}
Unless stated otherwise, $\K$ is a field of characteristic $0$, $n$ is a positive integer, and all Lie algebras and vector spaces are finite-dimensional over $\K$. An $n$-dimensional $\K$-vector space is denoted by $\K^n$, and $V$ denotes a certain $\K$-vector space. We write $\GL_n(\K)$ for the group of invertible $n \times n$ matrices and $\mathfrak{gl}_n(\K) \equiv \M_n(\K)$ for the Lie algebra of $n \times n$ matrices over $\K$. The symbol $U \oplus V$ denotes the direct sum of vector spaces, $\A \LO \B$ is the direct sum of Lie algebras $\A$ and $\B$, whereas $\A \oplus_\rho \B$ is the semidirect sum determined by an action $\rho \colon \B \to \Der(\A)$. Finally, $\s\{x_1, \ldots, x_n\}$ indicates the $n$-dimensional vector space with basis $(x_1, \ldots, x_n)$.
 
%----------Section 2-----------
\section{Preliminaries}\label{sec2}
%{Extensions of Lie algebras and trivial or semi-trivial splittings}

In this section, we recall extensions of Lie algebras and a special case, namely central extensions. Additionally, the notions of \emph{trivial} and \emph{semi-trivial} splittings are clarified.
%In this section, we first recall extensions of Lie algebras for later use. 

%------Subsection 2.1------
\subsection{Extensions of Lie algebras}

%--------Definition 2.1
\begin{definition}[Lie algebra extensions]
    Let $\A$, $\B$ and $\Cc$ be Lie algebras.
	\begin{enumerate}
		\item A sequence of Lie algebra homomorphisms of the form
		\begin{equation}%\label{E1}
		\xymatrix
	{
 		0 \ar[r] & \A \ar[r]^\iota & \B \ar[r]^\pi & \Cc \ar[r] & 0
 	}\tag{LE}          
 		\end{equation}
		is called a \emph{short exact sequence} of Lie algebras if $\iota$ is injective, $\pi$ is surjective and $\im \iota = \ker \pi$ is an ideal in $\B$. In this case, the sequence~\eqref{LE0} as well as the Lie algebra $\B$ is called a \emph{Lie algebra extension} (hereafter, an \emph{extension}) of $\Cc$ by $\A$.
		\item The following extension
		\begin{equation}\label{TE1}
        \xymatrix
        {
            0 \ar[r] & \A \, \ar@{^{(}->}[r]^{\hspace{-10pt}e} & \A \LO \Cc \ar@{->>}[r]^{\hspace{10pt}pr} & \Cc \ar[r] & 0
        }\tag{TLE}
		\end{equation}
		is called the \emph{trivial Lie algebra extension} (hereafter, \emph{trivial extension}) of $\Cc$ by $\A$. Here, $e \colon \A \hookrightarrow \A \LO \Cc$, $a \mapsto (a, 0)$ is the canonical embedding, and $pr \colon \A \LO \Cc \twoheadrightarrow \Cc$, $(a, c) \mapsto c$ is the canonical projection.         
	\end{enumerate}
\end{definition}

%------Definition 2.2----------
\begin{definition}[Equivalence of extensions]\label{DEE}
	Extension~\eqref{LE0} and the following extension
	\[
	\xymatrix
	{
 		0 \ar[r] & \A \ar[r]^{\bar{\iota}} & \bar{\B} \ar[r]^{\bar{\pi}} & \Cc \ar[r] & 0
 	}
	\]
	are said to be \emph{equivalent} if there exists a Lie algebra isomorphism $f \colon \B \to \bar{\B}$ making the following diagram commute
	\[
    \xymatrix
		{
			0 \ar[r] & \A \ar[r]^\iota \ar@{=}[d] & \B \ar[r]^\pi \ar[d]_{\raisebox{1.2ex}{$\scriptstyle\cong$}}^{\raisebox{1.5ex}{$\scriptstyle f$}} & \Cc \ar[r] \ar@{=}[d] & 0 \\			
			0 \ar[r] & \A \ar[r]^{\bar{\iota}} & \bar{\B} \ar[r]^{\bar{\pi}} & \Cc \ar[r] & 0.
		}
    \]
    From now on, every extension~\eqref{LE0} that is equivalent to a trivial extension~\eqref{TE1} will also be called a \emph{trivial extension} of $\Cc$ by $\A$.
\end{definition}

%-------------- Definition 2.3 ----------
\begin{definition}[Lie modules]\label{D113}
	Let $\A$ and $\B$ be two Lie algebras. If there is a Lie algebra homomorphism $\rho \colon \A \to \Der (\B)$ then $\B$ becomes an \emph{$\A$-module} by the structure $\A \times \B \to \B$ as follows  
    \[
       (a, x) \mapsto a \cdot x \co \rho(a) (x), \quad \text{for all } a \in \A \text{ and } x \in \B.
    \]
    If $\rho = 0$ then $\B$ is a \emph{trivial $\A$-module}, i.e., $a \cdot x = 0$ for all $a \in \A$ and $x \in \B$.
\end{definition}

%---------- Remark 2.4------------
\begin{remark}\label{R114}
For extension~\eqref{LE0}, we have the following observations. 
    \begin{enumerate}
        \item If we only consider extension~\eqref{LE0} in the category of vector spaces, then obviously $\B \cong \A \oplus \Cc$. Without loss of generality, we can consider $\B \equiv \A \oplus \Cc$ and each element $b \in \B$ can be written uniquely in the form $(a, c) \in \A \oplus \Cc$. Now $\iota \colon \A \hookrightarrow \B$ and $\pi \colon \B \twoheadrightarrow \Cc$ can be viewed as the canonical embedding and projection as follows
        \[
            \begin{array}{l l l l}
              \iota(a) \co (a, 0) & \text{and} & \pi(a, c) \co c, & \text{for all } a \in \A \text{ and } c\in \Cc. 
            \end{array}
        \]
        Thus, the linear embedding $\iota$ always has a left linear inverse $r \colon \B \to \A$, and the linear projection $\pi$ always has a right linear inverse $q \colon \Cc \to \B$ with
        \[
           \begin{array}{l l l l}
              r(a, c) \co a & \text{and} & q(c) \co (0, c), & \text{for all } a \in \A \text{ and } c\in \Cc.
            \end{array}
        \]
        
        \item\label{R114-2} Now we consider extension~\eqref{LE0} in the category of Lie algebras. If $r$ is also the left inverse of $\iota$ as a Lie algebra homomorphism then $r$ induces a trivial $\Cc$-module structure $\Cc \times \A \to \A$ on $\A$ as follows
        \[
            c \cdot a \co [r(0, c), a]_\A = [0, a]_\A = 0, \quad \text{for all } c \in \Cc \text{ and } a \in \A.
        \]
        Hence, extension~\eqref{LE0} becomes a trivial extension of $\Cc$ by $\A$, i.e., $\B \cong \A \LO \Cc$ (see~\cite{Che-Eil}).
        
        \item\label{R114-3} We continue to consider extension~\eqref{LE0} in the category of Lie algebras. If $q$ is also the right inverse of $\pi$ as a Lie algebra homomorphism then it also induces a $\Cc$-module structure on $\A$, but this structure is generally not trivial. Specifically, for all $a \in \A$ and $c \in \Cc$, we have 
        \[ \Cc \times \A \to \A; \, 
            (c, a) \mapsto c \cdot a \co \iota^{-1}([q(c), \iota(a)]_\B) = \iota^{-1}([(0,c), (a,0)]_\B).
        \]
        In other words, we get the action $\rho_q \colon \Cc \to \Der(\A)$ of $\Cc$ on $\A$ as follows
        \[
            \rho_q(c)(a) \co c \cdot a = \iota^{-1}([q(c), \iota(a)]_\B); \, \forall a \in \A, \, \forall c \in \Cc
        \]
        and $\B = \A \oplus_{\rho_q} \Cc$ via the action $\rho_q$ (see~\cite{Che-Eil}). That means $\B = \A \oplus \Cc$ and for all $(a_1, c_1), (a_2, c_2) \in \B$, their Lie bracket is defined as follows
        \[
            [(a_1, c_1), (a_2, c_2)]_{\B} := \Bigl([a_1, a_2]_{\A} + \rho_{q} (c_1)(a_2) - \rho_{q}(c_2)(a_1), \, [c_1, c_2]_{\mathcal{C}}\Bigr).
        \]
        Note that if $\iota$ maps $\A$ into the center of $\B$, i.e., $\im \iota = \ker \pi \subseteq \Z(\B)$, the action $\rho _q$ induced by $q$ will be trivial and extension~\eqref{LE0} also becomes a trivial extension of $\Cc$ by $\A$, i.e., $\B \equiv \A \LO \Cc$.% as the direct sum of Lie algebras.
    \end{enumerate}
\end{remark}

%-------Proposition 2.5------------
\begin{proposition}[{see \cite{Che-Eil}}]\label{P121}
	For an extension~\eqref{LE0}, the following statements are equivalent.
	\begin{enumerate}
		\item[(1)]\label{P25-1} It is a trivial extension of $\Cc$ by $\A$. 
		\item[(2)]\label{P25-2} The Lie algebra monomorphism $\iota$ has a left inverse $r \colon \B \to \A$.        
        \item[(3)]\label{P25-3} The Lie algebra epimorphism $\pi$ has a right inverse $q \colon \Cc \to \B$ and the action $\rho _q$ of $\Cc$ on $\A$ induced by $q$ is trivial.      
	\end{enumerate}
\end{proposition}

%------Subsection 2.2------
\subsection{Trivial and semi-trivial splittings of Lie algebra extensions}

%`In order to prove Theorem~\ref{thm1}, we need to present a precise definition of splitting and trivial splitting. 
Unlike the category of modules, in the category of Lie algebras, there is a distinction between trivial splitting and semi-trivial splitting of extensions. Specifically, we have the following definition.
%-------------Definition 2.6
\begin{definition}[Trivial or semi-trivial splittings of Lie algebra extensions]\label{DSE}
    We say that:    
    \begin{enumerate}
        \item Extension~\eqref{LE0} \emph{splits semi-trivially} (or simply, \emph{splits}) if the Lie algebra epimorphism $\pi$ has a right inverse $q$. In particular, then $\B = \A \oplus_{\rho_q} \Cc$ via the action $\rho _q$ of $\Cc$ on $\A$ induced by $q$.
        
        \item Extension~\eqref{LE0} \emph{splits trivially} if it satisfies one (and thus all) of the equivalence conditions stated in Proposition~\ref{P121}. By the first statement, when \eqref{LE0} splits trivially, we are also entitled to say that it is a {\it trivial extension} of $\Cc$ by $\A$.
    \end{enumerate}
\end{definition}

%----------Remark 2.7 ----------
\begin{remark}\label{R123}
    For the terms ``split trivially'' and ``split'' (semi-trivially), we have a few comments below.
    \begin{enumerate}
        \item Schottenloher~\cite[Definition 4.2]{Sch08} uses the term ``split'' which is actually ``split semi-trivially'' for every extension. But the word ``split'' for any central extension is ``split trivially'' as in Definition~\ref{DSE}.
        \item Weibel~\cite[Exercise 7.6.1]{Weibel} uses the term ``split'' which is actually ``split semi-trivially'' for every abelian Lie algebra extension as in Definition~\ref{DSE}.
        \item The terms ``trivial'' and ``split'' in Kerf et al.~\cite[Definition 18.1.6]{KBK97} correspond to ``split trivially'' and ``split semi-trivially'' in Definition~\ref{DSE}, respectively.
    \end{enumerate} 
\end{remark}

%---Subsection 2.3: Central Lie algebra extension---
%\subsection{Central Lie algebra extensions}

Next, we recall central extensions %, a class that naturally constitutes a restricted category of extensions in Lie theory
in the category of Lie algebras. They will be used in Section~\ref{sec4} to introduce the concept of cocomplete Lie algebras.
%-------------Definition 2.8----
\begin{definition}[{Central extensions~\cite[Definition 4.1]{Sch08}}]\label{DCLE}
A \emph{central Lie algebra extension} (for short, a \emph{central extension}) of $\Cc$ by $V$ is an extension
    \begin{equation}\label{CE1}
    \xymatrix
	{
 		0 \ar[r] & V \ar[r]^\iota & \B \ar[r]^\pi & \Cc \ar[r] & 0,
 	}\tag{CE}
\end{equation}
in which $\iota (V) \co \im \iota = \ker \pi$ is contained in the center $\Z(\B)$ of $\B$.
\end{definition}

%----Remark 2.9----
\begin{remark}\label{R29}
    Note that each central extension~\eqref{CE1} is a special case of generalized Lie algebra extensions~\eqref{LE0}. Therefore, one can also speak of the trivial and semi-trivial splittings of a central extension~\eqref{CE1}. However, by Item~\ref{R114-3} of Remark~\ref{R114}, these two notions of splitting coincide; in other words, every central extension that splits semi-trivially also splits trivially.
\end{remark}

%===================================================
%  SECTION 3
%===================================================
\section{Injective-type properties of complete Lie algebras}\label{sec3}

This section is devoted to complete Lie algebras and their injective-type behaviour in the category of Lie algebras. We first recall the definition and basic properties of complete Lie algebras, including their stability under finite direct sums. We then prove Theorem~\ref{thm1}, and discuss its relationship to the category-theoretic notion of an injective object: completeness yields a genuine lifting property under structural hypotheses, but is strictly weaker than the categorical notion (see Example~\ref{Ex-counter-cat-inj}). The section concludes with a practical algorithm for verifying completeness.

%-----Subsection 3.1 ---------
\subsection{Complete Lie algebras and their basic properties}
%In this subsection, we recall the concepts of {\it complete Lie algebras} that will serve as the foundation for the first main results of this paper. 

For a Lie algebra $\A$, we denote by $\ad(\A)$ the space of its inner derivations. Hence, $\Der(\A) \setminus\ad(\A)$ is precisely the set of all outer derivations of $\A$.

%---Definition 3.1 (Old 2.9) (Complete Lie Algebras)
\begin{definition}[{Complete Lie algebras~\cite[Chapter I]{Jac62}}]\label{D34} 
	A Lie algebra $\A$ is called \emph{complete} if $\Z(\A) = 0$ and $\Der(\A) = \ad (\A)$.
\end{definition}

%-------Remark 3.2 (Old 2.10)-----
\begin{remark}\label{R35}
    Complete Lie algebras have the following properties.
    \begin{enumerate}
        \item The trivial Lie algebra is complete.
        \item Every nonzero complete Lie algebra is non-abelian.
        \item All semisimple Lie algebras are always complete.
        \item\label{R35-4} All nonzero nilpotent Lie algebras are not complete.
    \end{enumerate}    
\end{remark}

%---Subsection 3.2: Complete Lie algebras and finite direct sums---
%\subsection{Complete Lie algebras and finite direct sums}

Unlike injective modules, which are not in general closed under direct sums, the class of complete Lie algebras is stable under finite direct sums. This stability will play a key role in the proof of Theorem~\ref{thm1}.

%----Theorem 3.5----Proposition 3.5---
\begin{proposition}[The stability of completeness under finite direct sums]\label{P36}
    Let $\A_1, \ldots, \A_m$ be Lie algebras with $m \geq 2$. Then $\A_1 \LO \cdots \LO \A_m$ is complete if and only if each $\A_i$ is complete.
\end{proposition}

\begin{proof}
    By finite induction, it suffices to establish the case $m = 2$. For convenience, we set $\A \co \A_1 \LO \A_2$.

    \begin{itemize}
        \item [$(\Leftarrow)$] If $\A_1$ and $\A_2$ are complete then $\Z(\A) = \Z(\A_1) \oplus \Z(\A_2) = \{0\}$. Furthermore, if $D \in \Der(\A)$, then $D(\A_i) \subseteq \A_i$ for $i = 1, 2$, since $\Z(\A_1) = \Z(\A_2) = 0$ and $[\A_1, \A_2]_\A = 0$; hence the restrictions $D|_{\A_1} \in \Der(\A_1)$ and $D|_{\A_2} \in \Der(\A_2)$ are well-defined. By completeness of each summand, there exist $a_1 \in \A_1$ and $a_2 \in \A_2$ such that
    $D|_{\A_1} = \ad_{a_1}$ and $D|_{\A_2} = \ad_{a_2}$. Setting
    $a \co (a_1, a_2) \in \A$, for every $x = (x_1, x_2) \in \A$,
    \[
        D(x) = \bigl(\ad_{a_1}(x_1), \ad_{a_2}(x_2)\bigr) = \ad_a(x),
    \]
    where the second equality uses $[\A_1, \A_2]_\A = \{0\}$. Hence
    $D = \ad_a \in \ad(\A)$, so $\Der(\A) = \ad(\A)$. Combined with
    $\Z(\A) = \{0\}$, this shows $\A$ is complete.

        \item [$(\Rightarrow)$] If $\A = \A_1 \LO \A_2$ is complete then the identity $\Z(\A) = \Z(\A_1) \oplus \Z(\A_2) = \{0\}$ forces $\Z(\A_i) = \{0\}$ for $i = 1, 2$. It remains to show that every derivation of $\A_i$ is inner. We argue for $\A_1$, the case of $\A_2$ is symmetric.
    
    Let $D_1 \in \Der(\A_1)$. Define $\tilde{D} \colon \A \to \A$ by
    \[
        \tilde{D}(x_1, x_2) \co (D_1(x_1), 0),
        \quad (x_1, x_2) \in \A_1 \LO \A_2.
    \]
    We claim that $\tilde{D} \in \Der (\A)$. In fact, for $x = (x_1, x_2), y = (y_1, y_2) \in \A$, using
    $[\A_1, \A_2]_\A = \{0\}$, one has
    \begin{align*}
        \tilde{D}([x, y]_\A)
         &= \tilde{D}\bigl([x_1, y_1]_{\A_1},\, [x_2, y_2]_{\A_2}\bigr)
          = \bigl(D_1([x_1, y_1]_{\A_1}),\, 0\bigr) \\
         &= \bigl([D_1(x_1), y_1]_{\A_1} + [x_1, D_1(y_1)]_{\A_1},\, 0\bigr) \\
         &= [\tilde{D}(x), y]_\A + [x, \tilde{D}(y)]_\A.
    \end{align*}
    Since $\A$ is complete, $\tilde{D} = \ad_a$ for some
    $a = (a_1, a_2) \in \A$. For every $x_1 \in \A_1$,
    \[
        (D_1(x_1), 0) = \tilde{D}(x_1, 0) = \ad_a(x_1, 0)
         = \bigl([a_1, x_1]_{\A_1},\, [a_2, 0]_{\A_2}\bigr)
         = (\ad_{a_1}(x_1),\, 0),
    \]
    where we used $[\A_2, \A_1]_\A = \{0\}$. Hence $D_1 = \ad_{a_1} \in \ad (\A_1)$, so $\Der(\A_1) = \ad(\A_1)$. Combined with $\Z(\A_1) = \{0\}$,
    this shows $\A_1$ is complete. By symmetry, $\A_2$ is also complete.
    \end{itemize}
    The proof of Proposition~\ref{P36} is complete.
\end{proof}

%----Subsection 3.3------
\subsection{Injective-type properties of complete Lie algebras}

We now prove the first main result of the paper, which concerns the injective-type properties of complete Lie algebras.

%----Proof of Theorem 1-----
\begin{proof}[{\bf Proof of Theorem~\ref{thm1}}]

We prove that $\eqref{Thm1-1} \Rightarrow \eqref{Thm1-2} \Rightarrow \eqref{Thm1-3} \Rightarrow \eqref{Thm1-1}$, which implies the equivalence of Items \eqref{Thm1-1}, \eqref{Thm1-2} and \eqref{Thm1-3}. The equivalence of Items \eqref{Thm1-3} and \eqref{Thm1-4} is established separately as it provides additional structural insight.

\medskip
%----Step 1: Prove (1) $\Rightarrow$ (2)----
\noindent\textbf{Step 1: Prove $\eqref{Thm1-1} \Rightarrow \eqref{Thm1-2}$.}
Assume $\A$ is complete and suppose $f, g$ as in the diagram \eqref{Diag1} and satisfy the assumptions of \eqref{Thm1-2}. Then, we can write $\A = \im g \LO J$ for some ideal $J$ of $\A$. Since $\A = \im g \LO J$ is complete, $\im g$ is complete by Proposition~\ref{P36}, i.e., $\Z(\im g) = \{0\}$ and $\Der(\im g) = \ad(\im g)$.
\begin{itemize}
    \item {\it The construction of the mapping} $h$. Since $\im f$ is an ideal of $\B$ by (\ref{Thm1-2i}), for each fixed $b \in \B$, the adjoint map $\ad_b$ restricts to an endomorphism of $\im f$. Transporting via the injection $f$, we obtain a linear map
\begin{equation}\label{eq:Db-def}
    D_b \colon \A' \to \A', \quad D_b(a') := f^{-1}\bigl([b, f(a')]_\B\bigr), \quad a' \in \A'.
\end{equation}
A direct check using the Jacobi identity in $\B$ shows that $D_b$ is a derivation of $\A'$. 
By \eqref{Thm1-2ii}, for every $b \in \B$ and every $a' \in \ker g$ we have 
\[\big([b, f(a')]_\B \in f(\ker g)\big) \Longrightarrow \big(D_b(\ker g) \subseteq \ker g\big).\]
Thus $D_b$ descends to a derivation of the quotient 
$\A'/\ker g \cong \im g$. We denote the induced derivation by
\[
    \bar{D}_b \in \Der\bigl(\A' / \ker g\bigr) \cong \Der(\im g) = \ad(\im g),
\]
and by $\ad^{\im g} \co \ad|_{\im g} \colon \im g \to \Der(\im g); \, a \, \mapsto \ad_{a}^{\im g} \in \Der(\im g)$ the adjoint action within the Lie algebra $\im g$. So for each $b \in \B$, there exists $a_b \in \im g$ with $\bar{D}_b = \ad_{a_b}^{\,\im g}$. 
Moreover, since $\Z(\im g) = \{0\}$, the element $a_b$ is uniquely determined. Define
\begin{equation}\label{eq:h-def}
    h \colon \B \to \A, \quad b \mapsto h(b) \co a_b \in \im g \subseteq \A.
\end{equation}
By construction, $\im h \subseteq \im g$.
    \item {\it Check the linearity of} $h$. The assignment $b \mapsto D_b$ is $\mathbb{K}$-linear by \eqref{eq:Db-def}, hence so is $b \mapsto \bar{D}_b$. Since $\ad^{\im g} \colon \im g \to \Der(\im g)$ is injective (as $\Z(\im g) = \{0\}$), the equality $\bar{D}_{b_1 + b_2} = \bar{D}_{b_1} + \bar{D}_{b_2} = \ad_{a_{b_1}}^{\im g} + \ad_{a_{b_2}}^{\im g} = \ad_{a_{b_1} + a_{b_2}}^{\im g}$ forces $a_{b_1 + b_2} = a_{b_1} + a_{b_2} \iff h(b_1 + b_2) = h(b_1) + h(b_2)$ for all $b_1, b_2 \in \B$. Similarly, it is also easy to see that $h(\alpha b) = \alpha h(b)$ for all $\alpha \in \K$ and $b \in \B$. 
    So $h$ is $\mathbb{K}$-linear.
    \item {\it Verify the equality} $h \circ f = g$. For any $a_0' \in \A'$ and any $a' \in \A'$,
\[
    D_{f(a_0')}(a') = f^{-1}\bigl([f(a_0'), f(a')]_\B\bigr) = [a_0', a']_{\A'}.
\]
Passing to the quotient and applying $g$, for every $a' \in \A'$,
\[
    \bar{D}_{f(a_0')}\bigl(g(a')\bigr) = g\bigl([a_0', a']_{\A'}\bigr) = [g(a_0'), g(a')]_\A = \ad_{g(a_0')}^{\im g}\bigl(g(a')\bigr),
\]
where the last equality uses the fact that $\im g$ is a Lie subalgebra of $\A$ even though it is actually an ideal of $\A$ by \eqref{Thm1-2ii}. Since $g$ corestricts to a surjection $\A' \twoheadrightarrow \im g$, this shows $\bar{D}_{f(a_0')} = \ad_{g(a_0')}^{\im g}$ on all of $\im g$. By the uniqueness of $a_b$ (note $g(a_0') \in \im g$), we conclude $h(f(a_0')) = g(a_0')$, i.e., $h \circ f = g$. This also implies $\im g \subset \im h$, hence, $\im h = \im g$.
    \item {\it Check the preservation of the Lie bracket of} $h$. The map $b \mapsto D_b$ is a Lie algebra homomorphism from $\B$ to $\Der(\A')$. Indeed, for $b_1, b_2 \in \B$ and $a' \in \A'$, the Jacobi identity in $\B$ gives
\begin{align*}
    D_{[b_1, b_2]_\B}(a') &= f^{-1}\bigl([[b_1, b_2]_\B, f(a')]_\B\bigr) \\
    &= f^{-1}\bigl([b_1, [b_2, f(a')]_\B]_\B - [b_2, [b_1, f(a')]_\B]_\B\bigr) \\
    &= D_{b_1}(D_{b_2}(a')) - D_{b_2}(D_{b_1}(a')) = [D_{b_1}, D_{b_2}]_{\Der(\A')}(a').
\end{align*}
Passing to the quotient $\A'/\ker g \cong \im g$, this yields 
\[\bar{D}_{[b_1, b_2]_\B} = [\bar{D}_{b_1}, \bar{D}_{b_2}]_{\Der(\im g)}.\]
Combining with $\bar{D}_b = \ad_{h(b)}^{\im g}$ and the fact that $\ad^{\im g}$ is a Lie algebra homomorphism, we obtain
\[
    \ad_{h([b_1, b_2]_\B)}^{\im g} = [\ad_{h(b_1)}^{\im g}, \ad_{h(b_2)}^{\im g}]_{\Der(\im g)} = \ad_{[h(b_1), h(b_2)]_{\im g}}^{\im g}.
\]
Since $\ad^{\im g}$ is injective, $h([b_1, b_2]_\B) = [h(b_1), h(b_2)]_{\im g} = [h(b_1), h(b_2)]_\A$, where the last equality uses that $\im g$ is a Lie subalgebra of $\A$. Hence, $h$ is a Lie algebra homomorphism.
\end{itemize}

%----Step 2: Prove (2) $\Rightarrow$ (3)----
\noindent\textbf{Step 2: Prove $\eqref{Thm1-2} \Rightarrow \eqref{Thm1-3}$.}
Let \eqref{LE0} be an arbitrary extension of $\Cc$ by $\A$. Apply \eqref{Thm1-2} to the $\A' \co \A$, $f \co \iota$ and $g \co \id_\A$. Then, \eqref{Thm1-2i}, \eqref{Thm1-2ii} and \eqref{Thm1-2iii} %of (\ref{Thm1-2}) 
are automatically satisfied since $\iota(\A) = \ker\pi$ is an ideal of $\B$ by definition of an extension, $\im \id_\A = \A$, which is trivially a direct Lie summand of itself (take $J = \{0\}$) and $\ker(\id_\A) = \{0\}$, so $\iota(\ker \id_\A) = \{0\}$ is the trivial ideal of $\B$. 
Hence, there exists a Lie algebra homomorphism $h \colon \B \to \A$ with $h \circ \iota = \id_\A$, i.e., $h$ is a left inverse of $\iota$. By Proposition~\ref{P121}, the extension \eqref{LE0} splits trivially.

\medskip
%----Step 3: Prove (3) $\Rightarrow$ (1)----
\noindent\textbf{Step 3: Prove $\eqref{Thm1-3} \Rightarrow \eqref{Thm1-1}$.}
Suppose every extension~\eqref{LE0} of any Lie algebra $\Cc$ by $\A$ splits trivially. We show the completeness of $\A$ by proving 
%\[\Der(\A) = \ad(\A); \quad \Z(\A) = \{0\}.\]
$\Der(\A) = \ad(\A)$ and $\Z(\A) = \{0\}$.
\begin{itemize}
    \item {\bf Part 3.1: Show $\Der(\A) = \ad(\A)$.} Let $D \in \Der(\A)$. Construct the semidirect sum $\B \co \A \oplus_D \mathbb{K}x$, where $x \notin \A$ and the Lie bracket extends that of $\A$ by
\[
    [x, a]_\B \co D(a), \quad \text{for all } a \in \A.
\]
This produces an extension
\[
    %0 \longrightarrow \A \longrightarrow \B \longrightarrow \mathbb{K}x \longrightarrow 0.
    \xymatrix
	{
 		0 \ar[r] & \A \ar[r] & \B \ar[r] & \K x \ar[r] & 0.
 	}
\]
By hypothesis, this extension splits trivially, so there exists an ideal $I$ of $\B$ with $\B = \A \LO I$. Since $\dim \B / \A = 1$, we have $\dim I = 1$ and write $I = \K y$.
Since $x \notin \A$ and $\B = \A \oplus I$ as a vector space, we may write $x = a_0 + \lambda y$ with $a_0 \in \A$ and $\lambda \in \mathbb{K} \setminus \{0\}$ (note $\lambda \neq 0$, otherwise $x \in \A$). Replacing $y$ by $\lambda y$ (which still spans $I$), we may assume $y = x - a_0$. Since $I$ is an ideal of $\B$ and $\A$ is also an ideal of $\B$, the intersection $I \cap \A = \{0\}$ forces $[y, \A]_\B \subseteq I \cap \A = \{0\}$. Thus for every $a \in \A$,
\[
    0 = [y, a]_\B = [x - a_0, a]_\B = [x, a]_\B - [a_0, a]_\A = D(a) - \ad_{a_0}(a),
\]
yielding $D = \ad_{a_0} \in \ad(\A)$. Hence $\Der(\A) = \ad(\A)$.
    \item {\bf Part 3.2: Show $\Z(\A) = \{0\}$.} We argue by contradiction. Suppose $\Z(\A) \neq \{0\}$ and pick $z \in \Z(\A) \setminus \{0\}$. Consider the vector space $\B := \A \oplus \mathbb{K} x \oplus \mathbb{K} y$ equipped with the bracket extending that of $\A$ by
\[
    [x, y]_\B = z, \quad [x, \A]_\B = [y, \A]_\B = 0.
\]
We first verify that $\B$ is a Lie algebra. Skew-symmetry is clear. For the Jacobi identity, consider the cases:
\begin{itemize}
    \item Triples in $\A$: Jacobi holds since $\A$ is a Lie algebra.
    \item Triples of the form $(x, y, a)$ with $a \in \A$:
    \[
        [[x, y]_\B, a]_\B + [[y, a]_\B, x]_\B + [[a, x]_\B, y]_\B = [z, a]_\A + 0 + 0 = 0,
    \]
    since $z \in \Z(\A)$.
    \item Triples of the form $(x, a_1, a_2)$ for $a_1, a_2 \in \A$:
    \[
        [[x, a_1]_\B, a_2]_\B + [[a_1, a_2]_\A, x]_\B + [[a_2, x]_\B, a_1]_\B = 0 + 0 + 0 = 0.
    \]
\end{itemize}
Similarly for $(y, a_1, a_2)$. So $\B$ is a Lie algebra. The inclusion $\A \hookrightarrow \B$ yields an extension
\[
    %0 \longrightarrow \A \longrightarrow \B \longrightarrow \B / \A \longrightarrow 0.
    \xymatrix
    {
    	0 \ar[r] & \A \ar[r] & \B \ar[r] & \B/\A \ar[r] & 0.
    }
\]
For the cosets $\bar{x} = x + \A$ and $\bar{y} = y + \A$ in $\B / \A$, we have $[\bar{x}, \bar{y}]_{\B/\A} = \overline{[x,y]}_\B = \bar{z} = 0$, so $\B / \A \cong \mathbb{K}^2$ as an abelian Lie algebra. By hypothesis, there exists an ideal $I \subseteq \B$ with $\B = \A \LO I$, where $I \cong \mathbb{K}^2$ is abelian. Let $x', y' \in I$ be the projections of $x, y$ onto $I$ along $\A$, so $x = x' + a_x$ and $y = y' + a_y$ with $a_x, a_y \in \A$. Since both $I$ and $\A$ are ideals with $I \cap \A = \{0\}$,
\[
    0 = [x', \A]_\B = [x - a_x, \A]_\B = [x, \A]_\B - [a_x, \A]_\A = -[a_x, \A]_\A,
\]
hence $a_x \in \Z(\A)$. Similarly $a_y \in \Z(\A)$. Since $I$ is abelian, $[x', y']_\B = 0$. Therefore,
\begin{align*}
    0 = [x', y']_\B &= [x - a_x, y - a_y]_\B \\
    &= [x, y]_\B - [x, a_y]_\B - [a_x, y]_\B + [a_x, a_y]_\A \\
    &= z - 0 - 0 + 0 = z,
\end{align*}
where the last identity uses $[x, a_y]_\B = [y, a_x]_\B = 0$ (by the construction of $\B$) and $[a_x, a_y]_\A = 0$ (since $a_x, a_y \in \Z(\A)$). This contradicts $z \neq 0$. Hence $\Z(\A) = \{0\}$.
\end{itemize}
\medskip
%----Step 4: Prove (3) $\Leftrightarrow$ (4)----
\noindent\textbf{Step 4: Prove $\eqref{Thm1-3} \Leftrightarrow \eqref{Thm1-4}$.}
\begin{itemize}
    \item {\bf Part 4.1: Show $\eqref{Thm1-3} \Rightarrow \eqref{Thm1-4}$.} Let $\iota \colon \A \hookrightarrow \B$ be a Lie algebra monomorphism with $\im\iota$ an ideal of $\B$. Setting $\Cc \co \B / \im\iota$ and letting $\pi \colon \B \twoheadrightarrow \Cc$ be the canonical projection, the sequence
\[
    %0 \longrightarrow \A \xrightarrow{\;\iota\;} \B \xrightarrow{\;\pi\;} \Cc \longrightarrow 0
    \xymatrix
	{
 		0 \ar[r] & \A \ar[r]^\iota & \B \ar[r]^\pi & \Cc \ar[r] & 0
 	}
\]
is an extension of $\Cc$ by $\A$. By \eqref{Thm1-3}, it splits trivially, so $\B = \im\iota \LO I$ for some ideal $I$.
    \item {\bf Part 4.2: Show $\eqref{Thm1-4} \Rightarrow \eqref{Thm1-3}$.} Given any extension
    \[
    %0 \longrightarrow \A \xrightarrow{\;\iota\;} \B \xrightarrow{\;\pi\;} \Cc \longrightarrow 0
    \xymatrix
	{
 		0 \ar[r] & \A \ar[r]^\iota & \B \ar[r]^\pi & \Cc \ar[r] & 0,
 	}
\]
    %$0 \to \A \xrightarrow{\iota} \B \xrightarrow{\pi} \Cc \to 0$, 
    the image $\im\iota = \ker\pi$ is an ideal of $\B$. By \eqref{Thm1-4}, $\im\iota$ is a direct summand: $\B = \im\iota \LO I$. Let $pr \colon \B \twoheadrightarrow I$ be the projection from $\B$ onto $I$. Then the composition $\pi \circ pr$ is an isomorphism $I \cong \Cc$, and one verifies directly that this realises the extension as a trivial one.
\end{itemize}
%The chain (\ref{Thm1-1}) $\Rightarrow$ (\ref{Thm1-2}) $\Rightarrow$ (\ref{Thm1-3}) $\Rightarrow$ (\ref{Thm1-1}) combined with (\ref{Thm1-3}) $\Leftrightarrow$ (\ref{Thm1-4})
Combining the above four steps completes the proof of Theorem \ref{thm1}.
\end{proof}

%-----Remark 3.3 about ``R-injective-type''-------
\begin{remark}\label{R-injective-type}
	As mentioned in Section~\ref{sec1}, Theorem~\ref{thm1} shows that complete Lie algebras act as “injective objects” for generalized extensions in the category of Lie algebras. Four equivalence conditions clarify this role from additional perspectives: condition~\eqref{Thm1-1} describes the internal structure ($A$ has a trivial center and only inner derivations); condition~\eqref{Thm1-2} describes the homomorphic lift property, extending from the familiarity of injective objects in the category of modules. Conditions~\eqref{Thm1-3} and \eqref{Thm1-4} are external separation descriptions; specifically, $A$ does not accept any nontrivial extensions as an ideal, which is equivalent to $A$ separating as an arbitrary direct component whenever it is an ideal. However, the full power of these equivalence properties is more subtle than the general analogy of injective objects in the category of modules. Specifically, condition~\eqref{Thm1-2} still requires the assumption that the image of $g$ must be a direct Lie summand of $\A$ -- a nontrivial structural constraint which, as Example~\ref{Ex-counter-cat-inj} shows, cannot be ignored. This also serves to delineate the precise gap between ``injective-type'' and full categorical injectivity in Lie algebra setting.  We will further clarify this in the next section.
\end{remark}

%-Example 3.4 (counterexample shows (2ii) very important- 
\begin{example}[A counterexample shows that \eqref{Thm1-2ii} cannot be ignored]\label{Ex-counter-cat-inj}
Let $\A = \mathfrak{sl}_2(\K) = \s \{e_1, e_2, e_3\}$ with
\[
    [e_1, e_2] = e_2, \quad
    [e_1, e_3] = -e_3, \quad
    [e_2, e_3] = e_1,
\]
as in Table~\ref{tab-completeLA}. Since $\A$ is semisimple, it is complete. Let $\B = \aff = \s \{x, y\}$ be the 2-dimensional affine Lie algebra with $[x, y] = y$. Set $\A' \co \K y \subset \B$ be a 1-dimensional ideal of $\B$. Since the inclusion $f \colon \A' \hookrightarrow \B$ is a Lie algebra monomorphism, assumption~\eqref{Thm1-2i} holds. Define a linear map
\[
    g \colon \A' \to \A, \quad g(y) \co e_1.
\]
Since $\dim \A' = 1$, the map $g$ is automatically a Lie algebra homomorphism, and $\ker g = \{0\}$, so assumption~\eqref{Thm1-2iii} that $f(\ker g)$ is an ideal of $\B$ is satisfied. However, $\im g = \K e_1$ is a proper 1-dimensional subspace of the simple Lie algebra $\mathfrak{sl}_2(\K)$. Since $\mathfrak{sl}_2(\K)$ does not admit any nontrivial Lie direct sum decomposition (its only ideals are $\{0\}$ and itself), $\im g$ is not a direct Lie summand of $\A$. Thus, assumption~\eqref{Thm1-2ii} %in (\ref{Thm1-2}) of Theorem \ref{thm1} 
fails.

We claim that no Lie algebra homomorphism $h \colon \B \to \A$ satisfies $h \circ f = g$. Assume for contradiction that such an $h$ exists, and set $X \co h(x) \in \A$. Applying $h$ to $[x, y] = y$ yields
\begin{equation}\label{eq-lift-constraint}
    [X, e_1] = [h(x), h(y)] = h([x, y]) = h(y) = e_1.
\end{equation}
On the other hand, writing $X = \alpha e_1 + \beta e_2 + \gamma e_3$ with $\alpha, \beta, \gamma \in \K$, a direct computation gives
\[
    [X, e_1] = -\beta e_2 + \gamma e_3 \in \s\{e_2, e_3\}.
\]
In particular, the $e_1$-component of $[X, e_1]$ vanishes for every $X \in \A$, contradicting~\eqref{eq-lift-constraint}. Hence, no such $h$ exists, even though
\begin{itemize}
    \item The target $\A$ is complete (in fact, semisimple);
    \item The monomorphism $f$ has image equal to an ideal of $\B$;
    \item The source $\A'$ is 1-dimensional abelian.
\end{itemize}
This confirms that the assumption~\eqref{Thm1-2ii} in Item~\eqref{Thm1-2} of Theorem~\ref{thm1} cannot be dropped, and consequently that categorical injectivity in the Lie algebra sense is strictly stronger than completeness.
\end{example}

%--Subsection 3.4: Complete Lie algebras and cohomological characterization---
%\subsection{Complete Lie algebras and cohomological characterization}
Now we turn to the cohomological viewpoint, which gives yet another equivalent formulation of completeness. Let $\ad \colon \A \to \mathfrak{gl}(\A)$ be the adjoint representation of $\A$. Then $\A$ becomes an $\A$-module under the adjoint action
$\A \times \A \to \A$, $(x, y) \mapsto x \cdot y \co \ad_x(y) = [x, y]$. Simple computations using the coboundary formula~\cite[Chapter~IV, Section~23, formula~(23.1)]{Che-Eil} yield
\[
    H^0(\A, \A) = \Z(\A) \quad \text{and} \quad
    H^1(\A, \A) = \Der(\A) / \ad(\A).
\]
From this point of view, complete Lie algebras admit an equivalent cohomological characterization as follows.
%---Proposition 3.6---
\begin{proposition}%[Cohomological characterization of complete Lie algebras]
\label{prop-H0andH1}
    A Lie algebra $\A$ is complete if and only if its first two cohomology groups with coefficients in itself vanish.%, i.e., $H^0(\A, \A) = H^1(\A, \A) = 0$.
\end{proposition}

%--Subsection 3.5: Identification of completeness---
\subsection{Identification of complete Lie algebras}

To the best of our knowledge, the classification of complete Lie algebras remains an open problem. Here we construct an algorithm for determining whether an $n$-dimensional Lie algebra $\A = \s\{e_1, \ldots, e_n\}$ is complete. Throughout this subsection, the structure constants $c_{ij}^k \in \K$ are defined by
\[
    [e_i, e_j] = \sum_{k=1}^{n} c_{ij}^k e_k, \quad 1 \leq i, j \leq n,
\]
with $c_{ii}^k = 0$ and $c_{ji}^k = -c_{ij}^k$.

For $x = \sum \limits_{j=1}^{n} \alpha_j e_j \in \A$, one has
\[
    [e_i, x] = \sum_{k=1}^{n} \biggl(\sum_{j=1}^{n} c_{ij}^k \alpha_j
              \biggr) e_k, \quad i = 1, \ldots, n.
\]
Therefore $x \in \Z(\A)$ if and only if
$(\alpha_1, \ldots, \alpha_n)$ is a solution of the homogeneous linear system
\begin{equation}\label{eq-center}
    \sum_{j=1}^{n} c_{ij}^k \alpha_j = 0, \quad
    1 \leq i, k \leq n.
\end{equation}
In particular, $\Z(\A) = \{0\}$ if and only if
system~\eqref{eq-center} admits only the trivial solution.

Now assume $\Z(\A) = \{0\}$. Since $\{\ad_{e_1}, \ldots, \ad_{e_n}\}$ generates $\ad(\A)$ and the map $x \mapsto \ad_x$ has kernel $\Z(\A) = \{0\}$, we have $\dim \ad(\A) = n$. Therefore $\A$ is complete if and only if $\Der(\A) = \ad(\A)$, i.e.,
\begin{equation}\label{der=ad}
    \dim \Der(\A) = n.
\end{equation}
Recall that a linear map $D \colon \A \to \A$ is a derivation of $\A$
if it satisfies the Leibniz rule
\[
    D([e_i, e_j]) = [D(e_i), e_j] + [e_i, D(e_j)], \quad
    1 \leq i < j \leq n.
\]
If $D(e_i) = \sum \limits_{p=1}^{n} z_{ip} e_p$ for $i = 1, \ldots, n$, i.e., the matrix of $D$ relative to the basis $(e_1, \ldots, e_n)$ is $[D] = (z_{ip})^T \in \M_n(\K)$, then expanding both sides of the Leibniz rule and comparing the coefficients of $e_p$ yields the
following equivalence: $D \in \Der(\A)$ if and only if
\begin{equation}\label{eq-derivation}
    \sum_{k=1}^{n} c_{ij}^k z_{kp}
     = \sum_{k=1}^{n} \bigl(c_{kj}^p\, z_{ik}
                            + c_{ik}^p\, z_{jk}\bigr),
    \quad 1 \leq i < j \leq n, \;\; 1 \leq p \leq n.
\end{equation}
System~\eqref{eq-derivation} is linear with unknowns $z_{ip}$, and $\dim \Der(\A)$ equals the number of its free parameters. In summary, we obtain Algorithm~\ref{alg1}.

%---Algorithm 1---
\begin{algorithm}[!h]
    \KwIn{Structure constants $c_{ij}^k \in \K$ of a Lie algebra $\A$}
    \KwOut{\texttt{True} if $\A$ is complete, otherwise \texttt{False}}
    \eIf{system~\eqref{eq-center} has nontrivial solutions}
        {return \texttt{False}}
        {Compute $\dim \Der(\A)$ via system~\eqref{eq-derivation}\;
         \eIf{equation~\eqref{der=ad} does not hold}
             {return \texttt{False}}
             {return \texttt{True}}
        }
    \caption{Identification of complete Lie algebras}\label{alg1}
\end{algorithm}

We illustrate Algorithm~\ref{alg1} with Example~\ref{Ex-alg1} below. For a given Lie algebra, we list only the nonzero Lie brackets.

%----Example 3.7----
\begin{example}\label{Ex-alg1}
Let $\A = \s\{e_1, e_2, e_3, e_4\}$ with
\[
    [e_1, e_3] = e_1, \quad
    [e_1, e_4] = -e_2, \quad
    [e_2, e_3] = e_2, \quad
    [e_2, e_4] = e_1.
\]
Here $c_{13}^1 = -c_{14}^2 = c_{23}^2 = c_{24}^1 = 1$, and
system~\eqref{eq-center} reduces to
$\alpha_1 = \alpha_2 = \alpha_3 = \alpha_4 = 0$, so
$\Z(\A) = \{0\}$. For $D \in \Der(\A)$ with matrix
$[D] = (z_{ip})^T \in \M_4(\K)$, system~\eqref{eq-derivation} becomes
\[
    \begin{cases}
        z_{13} = z_{14} = z_{23} = z_{24} = z_{33} = z_{34} = 0,\\
        z_{13} + z_{24} = 0,\\
        z_{14} - z_{23} = 0,\\
        z_{31} + z_{42} = 0,\\
        z_{32} - z_{41} = 0,\\
        z_{11} - z_{22} - z_{44} = 0,\\
        z_{11} - z_{22} + z_{44} = 0,\\
        z_{12} + z_{21} + z_{43} = 0,\\
        z_{12} + z_{21} - z_{43} = 0,
    \end{cases}
\]
which yields
\[
    [D] = \begin{bmatrix}
            a  &  b  & c  & d \\
           -b  &  a  & d  & -c \\
            0  &  0  & 0  & 0 \\
            0  &  0  & 0  & 0
          \end{bmatrix}, \quad a, b, c, d \in \K.
\]
Since $[D]$ contains four free parameters, $\dim \Der(\A) = 4$. By equation~\eqref{der=ad}, $\A$ is complete. Note that this is the Lie algebra $\s_{4,12}$ in Table~\ref{tab-completeLA}.
\end{example}

We have used Algorithm~\ref{alg1} to determine all complex and real complete Lie algebras up to dimension 4. The result is given in Table~\ref{tab-completeLA}.

%-----Table 1---
\begin{table}[!h]
    \centering
    \caption{Complex and real complete Lie algebras of dimension $\leq 4$}\label{tab-completeLA}
    \begin{tabular}{c p{7.7cm} c}
        \hline Lie algebras & \centering Non-zero Lie brackets & References \\
        \hline
            $\aff$ & $[e_1, e_2] = e_2$ & \\
        \hline
            $\mathfrak{sl}_2(\K)$ & $[e_1, e_2] = e_2$, $[e_1, e_3] = -e_3$, $[e_2, e_3] = e_1$ & \\
            $\mathfrak{so}_3(\R)$ & $[e_1, e_2] = e_3$, $[e_1, e_3] = -e_2$, $[e_2, e_3] = e_1$ & \\
        \hline
            $\aff \LO \aff$ & $[e_1, e_2] = e_2$, $[e_3, e_4] = e_4$ & \\
            $\Sf_{4,12}$ & $[e_1, e_3] = e_1$, $[e_1, e_4] = -e_2$, $[e_2, e_3] = e_2$, $[e_2, e_4] = e_1$ & \cite[Sec.~17.4]{SW14} \\
        \hline
    \end{tabular}
\end{table}

%====================================================
%  SECTION 4
%====================================================
\section{Cocomplete Lie algebras and projective-type properties with respect to central extensions}\label{sec4}

The investigation of objects with splitting properties in the category of Lie algebras leads to an interesting observation: although semisimple Lie algebras exhibit many features of projective objects, this category lacks true projective objects in the classical sense (Theorem~\ref{thm2}). To address this gap, we focus on the class of \emph{cocomplete} Lie algebras -- objects that play a ``projective-type'' role with respect to central extensions. We will characterize this class through the four equivalent conditions of Theorem~\ref{P-lifting-cocomplete} and, in the next subsections, by the vanishing of the second cohomology group with trivial coefficients (Proposition~\ref{P321}), the latter providing Algorithm~\ref{alg2} as a computational tool. Building on this foundation, we then present a detailed classification of the almost abelian case. The starting point is the following well-known consequence of Levi's
theorem~\cite[Part~I, Chapter~VI, Theorem~4.1]{Serre}.

%---Proposition 4.1 (Levi's theorem)
\begin{proposition}[Levi's theorem]\label{P124}
Let $\pi \colon \B \to \Cc$ be a surjective Lie algebra homomorphism onto a semisimple Lie algebra $\Cc$. Then there exists a Lie algebra homomorphism $q \colon \Cc \to \B$ such that $\pi \circ q = \id_\Cc$.
\end{proposition}

%---Subsection 4.1: Proof of Theorem 2 and cocomplete Lie algebras----
\subsection{Proof of Theorem~\ref{thm2} and cocomplete Lie algebras}\label{subsec:4-1}
By Levi's theorem, if $\Cc$ is semisimple then every
extension~\eqref{LE0} splits semi-trivially. However, if we replace ``splits semi-trivially'' by ``splits trivially,'' Levi's theorem is no longer valid: this is the content of Theorem~\ref{thm2}, whose proof we now give.

\begin{proof}[{\bf Proof of Theorem~\ref{thm2}}]
    Assume that $\Cc$ is a nontrivial Lie algebra. By Ado's theorem~\cite[Chapter~V, \S8, p.\,153]{Serre}, there exists a faithful representation
$\rho \colon \Cc \to \mathfrak{gl}(V) \equiv \Der(V)$ of $\Cc$ in some vector space $V$. Choose $\A \co V$ as an abelian Lie algebra. Since $\Cc \neq \{0\}$ and $\rho$ is faithful, there exist $a_0 \in \A$ and $c_0 \in \Cc$ such that $\rho(c_0)(a_0) \neq 0$. Then $\rho$ defines a (nontrivial) $\Cc$-module structure on $\A$
via
\[
    \Cc \times \A \to \A, \quad
    (c, a) \mapsto c \cdot a \co \rho(c)(a).
\]
Define the Lie algebra $\B \co \A \oplus_\rho \Cc$, i.e.,
$\B = \A \oplus \Cc$ as a vector space with bracket
\[
    [(a_1, c_1), (a_2, c_2)]_\B
     \co \bigl(c_1 \cdot a_2 - c_2 \cdot a_1,\,
              [c_1, c_2]_\Cc\bigr).
\]
The canonical embedding $\iota \colon \A \to \B$,
$a \mapsto (a, 0)$, and projection $\pi \colon \B \to \Cc$, $(a, c) \mapsto c$, fit into an extension~\eqref{LE0} of $\Cc$ by
$\A$. Since $[(0, c_0), (a_0, 0)]_\B = (c_0 \cdot a_0, 0) \neq 0$, we have $[\Cc, \A]_\B \neq \{0\}$, hence
$\B = \A \oplus_\rho \Cc \ncong \A \LO \Cc$. Therefore, this extension does not split trivially. The proof is complete.  
\end{proof}

The negative result of Theorem~\ref{thm2} confirms that a global ``projective-type'' property cannot exist in the category of Lie algebras. To pursue a meaningful projective analogue, we must restrict attention to the natural and extensively studied class of central extensions, as anticipated in Section~\ref{sec1}. Namely, we have the following definition.

%-Definition 4.2(Old 4.3):Cocomplete Lie algebras---
\begin{definition}[Cocomplete Lie algebras]\label{D311}
A Lie algebra $\Cc$ is called {\it cocomplete} if every central extension (\ref{CE1}) of $\Cc$ by any vector space $V$ splits trivially.
\end{definition}

The concept of the cocomplete Lie algebras provides a natural counterpart to complete Lie algebras within the framework of central extensions. Hence, the remainder of this section is devoted to a detailed study of this class of Lie algebras and its structural properties.

%---Subsection 4.2: Cocompleteness and quasi-projective property -- Proof of Theorem 3----
\subsection{Cocompleteness and projective-type property: Proof of Theorem~\ref{P-lifting-cocomplete}}\label{subsec:4-2}

We now turn to the proof of Theorem~\ref{P-lifting-cocomplete}, which establishes the equivalence of four formulations of the ``projective-type'' property: the splitting condition for central extensions, the direct-summand condition for central kernels, and a lifting property along central epimorphisms.

\begin{proof}[\bf Proof of Theorem~\ref{P-lifting-cocomplete}]

By Definition~\ref{D311}, the equivalence $\eqref{Thm3-0} \Leftrightarrow \eqref{Thm3-1}$ is immediate. It suffices to establish the two equivalences $\eqref{Thm3-1} \Leftrightarrow \eqref{Thm3-2}$ and $\eqref{Thm3-1} \Leftrightarrow \eqref{Thm3-3}$.
\vskip 0.15cm
%\medskip
%---Step 1: Prove (1) $\Leftrightarrow$ (2)--
\noindent\textbf{Step 1: Prove $\eqref{Thm3-1} \Leftrightarrow \eqref{Thm3-2}$.}
\begin{itemize}
    \item {\bf Show $\eqref{Thm3-1} \Rightarrow \eqref{Thm3-2}$}. Assume that \eqref{Thm3-1} in Theorem~\ref{P-lifting-cocomplete} holds. We need to prove that \eqref{Thm3-2} in Theorem~\ref{P-lifting-cocomplete} holds. Let $\pi \colon \B \twoheadrightarrow \Cc$ be a Lie algebra epimorphism with $\ker \pi \subseteq \Z(\B)$. We will show that $\ker \pi$ is a direct summand of $\B$. By setting $V \co \ker \pi$, which is abelian as a subspace of the center, the sequence
    \[
       \xymatrix
	{
 		0 \ar[r] & V \ar@{^{(}->}[r]^\iota & \B \ar[r]^\pi & \Cc \ar[r] & 0
 	}
    \]
    is a central extension of $\Cc$ by $V$. By \eqref{Thm3-1}, it splits trivially, so by Proposition~\ref{P121}, $V = \ker \pi$ is a direct summand of $\B$, i.e., $\B = \ker \pi \LO I$ for some ideal $I$ of $\B$, with $I \cong \Cc$.
    
    \item {\bf Show $\eqref{Thm3-2} \Rightarrow \eqref{Thm3-1}$.} Assume that \eqref{Thm3-2} in Theorem~\ref{P-lifting-cocomplete} holds. We need to show that \eqref{Thm3-1} in Theorem~\ref{P-lifting-cocomplete} holds. Given any central extension~\eqref{CE1}, the kernel $\ker \pi = \iota(V)$ is always contained in $\Z(\B)$. By \eqref{Thm3-2}, $\ker \pi$ is a direct summand of $\B$, i.e., $\B = \ker \pi \LO I \cong V \LO \Cc$. Hence, the extension splits trivially.
\end{itemize}

\smallskip
%---Step 2: Prove (1) $\Leftrightarrow$ (3)--
\noindent\textbf{Step 2: Prove $\eqref{Thm3-1} \Leftrightarrow \eqref{Thm3-3}$}.
\begin{itemize}
    \item {\bf Show $\eqref{Thm3-1} \Rightarrow \eqref{Thm3-3}$.} Assume \eqref{Thm3-1} in Theorem~\ref{P-lifting-cocomplete} holds. We need to prove that \eqref{Thm3-3} in Theorem~\ref{P-lifting-cocomplete} holds. Consider two Lie algebra homomorphisms $f, g$ as in a diagram of the form~\eqref{Diag2} satisfying \eqref{Thm3-3i} and \eqref{Thm3-3ii} in Theorem~\ref{P-lifting-cocomplete}. We need to show that there exists a Lie algebra homomorphism $h \colon \Cc \to \B$ such that $f \circ h = g$. Indeed, we set 
    \[
        E \co \bigl\{(b, c) \in \B \times \Cc \mid f(b) = g(c) \in \Cc'\bigr\} \subseteq \B \times \Cc.
    \]
    For any $(b_1, c_1), (b_2, c_2) \in E$, we have
    \[
        f([b_1, b_2]_\B) = [f(b_1), f(b_2)]_{\Cc'}
        = [g(c_1), g(c_2)]_{\Cc'} = g([c_1, c_2]_\Cc),
    \]
    so $([b_1, b_2]_\B, [c_1, c_2]_\Cc) \in E$. Hence, $E$ is a Lie subalgebra of $\B \times \Cc$.

    Consider the canonical projections 
\[pr_{\B} \colon \B \times \Cc \twoheadrightarrow \B, (b, c) \mapsto pr_{\B}(b, c) \co b,\]
\[pr_{\Cc} \colon \B \times \Cc \twoheadrightarrow \Cc, (b, c) \mapsto pr_{\Cc}(b, c) \co c.\]
%$pr: \B \times \Cc \twoheadrightarrow \Cc; \, (b, c) \mapsto pr(b, c) : = c$ and 
    Now we set $\pi_E \co pr_{\Cc}|_{E} \colon E \to \Cc$. %, $(c, b) \mapsto c$. 
    Since $f \colon \B \to \Cc'$ is surjective, for every $c \in \Cc$, there exists $b \in \B$ with $f(b) = g(c)$, in particular $(b, c) \in E$ and $\pi_E$ is also surjective. Moreover,
\[
    \ker \pi_E
     = \{(b, 0) \in E \mid f(b) = g(0) = 0\}
     = \ker f \times \{0\} \subseteq \Z(\B) \times \{0\},
\]
i.e., $\ker \pi_E$ is abelian. In particular, $\ker \pi_E \subseteq \Z(E)$. Hence,
    \begin{equation}\label{eq-pullback-ext}
    %0 \to \ker \pi_E \hookrightarrow E       \xrightarrow{\;\pi_E\;} \Cc \to 0
       \xymatrix
	{
 		0 \ar[r] & \ker \pi_E \ar@{^{(}->}[r]^{\hspace{8pt}\iota} & E \ar[r]^{\pi_E} & \Cc \ar[r] & 0
 	  }
    \end{equation}
    is a central extension of $\Cc$ by the abelian Lie algebra $\ker \pi_E$.

    By~\eqref{Thm3-1} in Theorem~\ref{P-lifting-cocomplete}, the central extension~\eqref{eq-pullback-ext} splits trivially. By Proposition~\ref{P121}, $\pi_E$ admits a Lie algebra right inverse $s \colon \Cc \to E$ with $\pi_E \circ s = \id_\Cc$. Set $h \co (pr_{\B}|_{E}) \circ s \colon \Cc \to \B$ as composition of two Lie algebra
    homomorphisms. We just need to check that $f \circ h = g$. Indeed, from $\pi_E \circ s = \id_\Cc$ we must have $s(c) = (b, c) \in E$ for some $b \in \B$ with $f(b) = g(c)$. Since $h \co (pr_{\B}|_{E}) \circ s$, we must have $h (c) = b$. Thus, $f(h(c)) = f(b) = g(c)$ for all $c \in \Cc$. Therefore, $f \circ h = g$.
    
    \item {\bf Show $\eqref{Thm3-3} \Rightarrow \eqref{Thm3-1}$.} Assume \eqref{Thm3-3} in Theorem~\ref{P-lifting-cocomplete} holds. We need to prove that \eqref{Thm3-1} in Theorem~\ref{P-lifting-cocomplete} holds.
    Let~\eqref{CE1} be an arbitrary central extension of $\Cc$ by an abelian Lie algebra $V$. Applying~\eqref{Thm3-3} in Theorem~\ref{P-lifting-cocomplete} with $f \co \pi$, $\Cc' \co \Cc$, and $g \co \id_\Cc$ (which is valid since $\ker \pi = \iota(V) \subseteq \Z(\B)$), we obtain a Lie algebra homomorphism $h \colon \Cc \to \B$ with $\pi \circ h = \id_\Cc$, i.e., $h$ is a right inverse of $\pi$. Since the extension is central, $\iota(V) \subseteq \Z(\B)$, so the induced action $\rho_h$ of $\Cc$ on $V$ is trivial by Item~\ref{R114-3} of Remark~\ref{R114}. By Proposition~\ref{P121}, the extension splits trivially.
\end{itemize}
    Combining the above two steps completes the proof of Theorem~\ref{P-lifting-cocomplete}.
\end{proof}

%---Remark 4.3----
\begin{remark}\label{R-asymmetry-inj-proj}
The pair of Theorems~\ref{thm1} and~\ref{P-lifting-cocomplete} reveals a fundamental asymmetry between the injective-type and projective-type properties in the category of Lie algebras. The property of injective-type lifting homomorphism \eqref{Thm1-2} in Theorem~\ref{thm1}~requires three structural conditions \eqref{Thm1-2i}, \eqref{Thm1-2ii} and \eqref{Thm1-2iii} on the
pair $(f, g)$; in particular, the image $\im g$ must be a direct Lie summand of the complete target $\A$. The projective-type lifting homomorphism \eqref{Thm3-3} in Theorem~\ref{P-lifting-cocomplete}, by contrast, imposes \emph{no} additional condition on the target homomorphism $g \colon \Cc \to \Cc'$; the centrality condition $\ker f \subseteq \Z(\B)$ on the epimorphism alone is sufficient. This asymmetry is intrinsic to the non-abelian structure of the category: the center and the derivation algebra of $\A$ govern the obstructions to extending homomorphisms \emph{into} $\A$, whereas centrality of the kernel alone controls the obstructions to lifting homomorphisms \emph{out of} $\Cc$.
\end{remark}

%-Subsection 4.3:Cocomplete Lie algebras and cohomology-
\subsection{Cohomological characterizations of cocomplete Lie algebras}\label{subsec:4-3}

%We now introduce the central notion of the section, justified by the equivalence of the conditions in Theorem~\ref{P-lifting-cocomplete}.

%---Remark 4.4---
%\begin{remark}\label{R-cocomplete-categorical}
In categorical language, condition~\eqref{Thm3-2} in Theorem~\ref{P-lifting-cocomplete} indicates that cocomplete Lie algebras are precisely the projective objects with respect to ``central'' epimorphisms in the category of Lie algebras. This is the ``quasi projective-type'' correspondence that naturally arises within the context of the limitation on central extensions and is somewhat the ``dual'' to the ``injective-type'' role assumed by complete Lie algebras through Theorem~\ref{thm1}.
%\end{remark}

%---Proposition 4.5----
\begin{proposition}\label{P312}
    Every semisimple Lie algebra is cocomplete.
\end{proposition}

\begin{proof}
    Let $\Cc$ be semisimple and \eqref{CE1} a central extension of $\Cc$ by $V$. By Proposition~\ref{P124}, the epimorphism $\pi$ admits a Lie algebra right inverse $q \colon \Cc \to \B$, so $\B = V \oplus_{\rho_q} \Cc$ for the action $\rho_q \colon \Cc \to \Der(V)$ induced by $q$. Since $V$ is abelian and $\iota(V) \subseteq \Z(\B)$, the action $\rho_q$ is trivial by Item~\ref{R114-3} of Remark~\ref{R114}. Hence, the central extension~\eqref{CE1} splits trivially, i.e., $\Cc$ is cocomplete.
\end{proof}

The next result provides a computationally necessary and sufficient condition for verifying whether a Lie algebra $\Cc$ is cocomplete. To formulate and prove it, we first need the following lemma.

%---Lemma 4.6----
\begin{lemma}\label{L-trivial-coeff}
For any Lie algebra $\Cc$ and every vector space $V$ regarded as a
trivial $\Cc$-module, we have
\[
    H^p(\Cc, V) \cong H^p(\Cc, \K) \otimes_\K V,
    \quad p \in \N.
\]
\end{lemma}

\begin{proof}
Since $V$ is a trivial $\Cc$-module,
$C^p(\Cc, V) \cong C^p(\Cc, \K) \otimes_\K V$. Let
\[d^p_V \colon C^p(\Cc, V) \to C^{p+1}(\Cc, V); \, d^p_\K \colon \Lambda^p(\Cc^*) \to \Lambda^{p+1}(\Cc^*)\]
%$d^p_V \colon C^p(\Cc, V) \to C^{p+1}(\Cc, V)$ and
%$d^p_\K \colon \Lambda^p(\Cc^*) \to \Lambda^{p+1}(\Cc^*)$ 
be the coboundary operators. Then $d^p_V = d^p_\K \otimes \id_V$, so
\[
    Z^p(\Cc, V) \cong \ker d^p_\K \otimes_\K V
    \quad \text{and} \quad
    B^p(\Cc, V) \cong \im d^{p-1}_\K \otimes_\K V \quad (p > 0).
\]
Therefore,
\[
    H^p(\Cc, V)
     \cong \frac{\ker d^p_\K \otimes_\K V}
                {\im d^{p-1}_\K \otimes_\K V}
     \cong \frac{\ker d^p_\K}{\im d^{p-1}_\K} \otimes_\K V
     = H^p(\Cc, \K) \otimes_\K V. \qedhere
\]
\end{proof}

%----Proposition 4.7----
\begin{proposition}\label{P321}
A Lie algebra $\Cc$ is cocomplete if and only if
$H^2(\Cc, \K) = 0$.
\end{proposition}

\begin{proof}
    By Schottenloher~\cite[Chapter~IV, Remark~4.7]{Sch08}, a central extension~\eqref{CE1} of $\Cc$ by $V$ splits trivially if and only if its cohomology class in $H^2(\Cc, V)$ vanishes. Hence, every central extension of $\Cc$ by $V$ splits trivially if and only if $H^2(\Cc, V) = 0$. By Lemma~\ref{L-trivial-coeff}, this holds for every vector space $V$ if and only if $H^2(\Cc, \K) = 0$.
\end{proof}

Recall that $b_p(\Cc) \co \dim H^p(\Cc, \K)$ is the
\emph{$p$-th Betti number} of $\Cc$.

%----Corollary 4.8----
\begin{corollary}\label{C-betti}
A Lie algebra $\Cc$ is cocomplete if and only if $b_2(\Cc) = 0$.
\end{corollary}

%-Subsection 4.4: Cocomplete Lie algebras and finite direct sums----
%\subsection{Cocomplete Lie algebras and finite direct sums}\label{subsec:4-31}
Recall that a Lie algebra $\Cc$ is \emph{perfect} if $\Cc = [\Cc, \Cc]$. The class of cocomplete Lie algebras exhibits a property analogous to that of projective modules under direct sums, although a suitable restriction is required. More precisely, the following result can be viewed as a ``dual'' version of Proposition~\ref{P36}.

%----Proposition 4.9 (4.6 Old)-----
\begin{proposition}\label{P43}
    Let $\Cc_1, \ldots, \Cc_m$ be Lie algebras with $m \geq 2$. Then $\Cc_1 \LO \cdots \LO \Cc_m$ is cocomplete if and only if the following two conditions hold.
    \begin{enumerate}[(1)]
        \item All $\Cc_1, \ldots, \Cc_m$ are cocomplete.
        \item At least $m-1$ of the summands $\Cc_1, \ldots, \Cc_m$ are perfect.
    \end{enumerate}
\end{proposition}

\begin{proof}
    By finite induction, it suffices to establish the case $m = 2$. By the K\"unneth formula~\cite[Exercise~7.3.8]{Weibel}, one has
    \[
        H^2 \bigg(\Cc_1 \LO \Cc_2, \K\bigg)
        \cong H^2(\Cc_1, \K) \oplus \big[ H^1(\Cc_1, \K) \otimes H^1(\Cc_2, \K)\big] \oplus H^2(\Cc_2, \K).
    \]
    By Proposition~\ref{P321}, $\Cc_1 \LO \Cc_2$ is cocomplete if and only if %$H^2 \bigg(\Cc_1 \LO \Cc_2, \K\bigg) = 0$, which is equivalent to
    \[
        H^2(\Cc_1, \K) = H^1(\Cc_1, \K) \otimes H^1(\Cc_2, \K) = H^2(\Cc_2, \K) = 0.
    \]
    The vanishing of the first and third components above means that $\Cc_1$ and $\Cc_2$ are cocomplete. Furthermore, the vanishing of the second component means $H^1(\Cc_1, \K) = 0$ or $H^1(\Cc_2, \K) = 0$. Since $H^1(\Cc, \K) \cong (\Cc/[\Cc, \Cc])^*$, we must have $\Cc_1 = [\Cc_1, \Cc_1]$ or $\Cc_2 = [\Cc_2, \Cc_2]$, i.e., $\Cc_1$ or $\Cc_2$ is perfect.
\end{proof}

%\noindent {\it Proof}. 
%\begin{proof}
%\begin{itemize}
%    \item Show ($\Longrightarrow$): Let $\Cc = \Cc_1 \LO \Cc_2$ be cocomplete. Then $H^2(\Cc, \K) = 0$. On the other hand, the K\"unneth formula~\cite[Exercise~7.3.8]{Weibel} gives
%    \[H^2(\Cc, \K)
%     \cong H^2(\Cc_1, \K)
%           \oplus \bigl(H^1(\Cc_1, \K) \otimes
%                        H^1(\Cc_2, \K)\bigr)
%           \oplus H^2(\Cc_2, \K).\]
%    Since $H^2(\Cc, \K) = 0$, it is mandatory to have
%    \begin{equation*}
%        \begin{cases}
%            H^2(\Cc_1, \K) = H^2(\Cc_2, \K) = 0, \\
%            H^1(\Cc_1, \K) = 0 \, \vee \, H^1(\Cc_2, \K) = 0.
%        \end{cases}
%    \end{equation*}
%    So all $\Cc_1, \Cc_2$ are cocomplete and at least one of $\Cc_1, \Cc_2$ is perfect. 
%    \item Show ($\Longleftarrow$): Assume that $\Cc_1, \Cc_2$ are two cocomplete Lie algebras, i.e. $H^2(\Cc_1, \K) =$ $H^2(\Cc_2, \K) = 0$. Furthermore, at least one of them is perfect. Without loss of generality, assume that $\Cc_2$ is perfect, i.e. $\Cc_2 = [\Cc_2, \Cc_2]$. Then
%    \[H^1(\Cc_2, \K) \cong (\Cc_2/[\Cc_2, \Cc_2])^* = 0 \Rightarrow H^1(\Cc_1, \K) \otimes
 %                       H^1(\Cc_2, \K) = 0.\]
%    From there, applying the K\"unneth formula~\cite[Exercise~7.3.8]{Weibel}, we get $H^2(\Cc_1 \LO \Cc_2, \K) = 0$. This also means that $\Cc = \Cc_1 \LO \Cc_2$ is a cocomplete Lie algebra. The proof is complete. $\hfill \square$   
%\end{itemize}

%-Subsection 4.5: Identification of cocomplete Lie algebras
\subsection{Identification of cocomplete Lie algebras}\label{subsec:4-32}
The classification of cocomplete Lie algebras is, in general, not straightforward and will be a subject of future work. First, regarding the ``magnitude'' of the class of cocomplete Lie algebras, we have the following observations.

%----Remark 4.10----
\begin{remark}\label{R-cocomplete-examples}
\leavevmode
\begin{itemize}
    \item The class of cocomplete Lie algebras is strictly larger than that of semisimple Lie algebras. For instance, the affine Lie algebra $\aff$ satisfies $H^2(\aff, \K) = 0$, hence is cocomplete, but is evidently not semisimple.

    \item According to Dixmier~\cite[Th\'eor\`eme~2]{Dix55}, if $\Cc$ is an $n$-dimensional nilpotent Lie algebra then
    \[
        b_0(\Cc) \geq 1, \quad b_n(\Cc) \geq 1, \quad
        b_p(\Cc) \geq 2 \text{ for all } 0 < p < n.
    \]
    It is also well-known that $b_2(\K^n) = \binom{n}{2}$. By Corollary~\ref{C-betti}, every nilpotent Lie algebra of dimension $n \geq 2$ fails to be cocomplete. 
    
    \item The direct sum $\K^n \LO \Cc$ is not  cocomplete for $n \geq 2$, regardless of $\Cc$.
\end{itemize}
\end{remark}

Now we construct an algorithm to identify whether an $n$-dimensional Lie algebra $\Cc = \s\{e_1, \ldots, e_n\}$ with structure constants $c_{ij}^k$ is cocomplete. Let $\Lambda^1(\Cc^*) \co \Cc^* = \s\{e_1^*, \ldots, e_n^*\}$ be the dual space of $\Cc$. The exterior differential (or coboundary operator) of $e_i^*$ is
\begin{equation}\label{exterior-derivative}
    d e_i^* = -\frac{1}{2}
                \sum_{j, k = 1}^{n}
                    c_{jk}^i\, e_j^* \wedge e_k^*.
\end{equation}
The space of 2-forms is $\Lambda^2(\Cc^*)
 = \s\{e_i^* \wedge e_j^* \mid 1 \leq i < j \leq n\}$, with 
 \begin{equation}\label{second-exterior-derivative}
    d(e_i^* \wedge e_j^*)
     = d e_i^* \wedge e_j^* - e_i^* \wedge d e_j^*.
\end{equation}
A 2-form $\omega \in \Lambda^2(\Cc^*)$ is \emph{closed} if $d \omega = 0$, and \emph{exact} if $\omega = d \theta$ for some
$\theta \in \Lambda^1(\Cc^*)$. Then
\begin{equation}\label{Z2(g)andB2(g)}
    \begin{cases}
        Z^2(\Cc, \K)
         = \{\omega \in \Lambda^2(\Cc^*)
              \mid \omega \text{ is closed}\}, \\
        B^2(\Cc, \K)
         = \{\omega \in \Lambda^2(\Cc^*)
              \mid \omega \text{ is exact}\}.
    \end{cases}
\end{equation}
By Proposition~\ref{P321}, $\Cc$ is cocomplete if and only if 
\[H^2(\Cc, \K) = 0 \iff Z^2(\Cc, \K) = B^2(\Cc, \K).\]
%$H^2(\Cc, \K) = 0$, i.e., $Z^2(\Cc, \K) = B^2(\Cc, \K)$. 
This forms the basis for Algorithm~\ref{alg2} below.

%---Algorithm 2----
\begin{algorithm}[h]
    \KwIn{Structure constants $c_{ij}^k \in \K$ of a Lie algebra
           $\Cc$}
    \KwOut{\texttt{True} if $\Cc$ is cocomplete, otherwise
           \texttt{False}}
    Compute $Z^2(\Cc, \K)$ and $B^2(\Cc, \K)$
    by~\eqref{Z2(g)andB2(g)}\;
    \eIf{$Z^2(\Cc, \K) = B^2(\Cc, \K)$}
        {return \texttt{True}}
        {return \texttt{False}}
    \caption{Identification of cocomplete Lie algebras}\label{alg2}
\end{algorithm}

%----Example 4.11-----
\begin{example}\label{Ex-cocomplete}
Let $\Cc = \s \{e_1, e_2, e_3, e_4\}$ with
\[
    [e_1, e_4] = 2 e_1,\; [e_2, e_4] = e_2,\; [e_3, e_4] = e_2 + e_3.
\]
Here $c_{14}^1 = 2$ and $c_{24}^2 = c_{34}^2 = c_{34}^3 = 1$, and
$\Lambda^1(\Cc^*) = \s\{e_1^*, e_2^*, e_3^*, e_4^*\}$. By
formula~\eqref{exterior-derivative},
\[
    d e_1^* = -2 e_1^* \wedge e_4^*, \;
    d e_2^* = -e_2^* \wedge e_4^* - e_3^* \wedge e_4^*, \;
    d e_3^* = -e_3^* \wedge e_4^*, \;
    d e_4^* = 0.
\]
For $\omega = \sum_{i = 1}^{4} a_i e_i^* \in \Lambda^1(\Cc^*)$,
one has
\[d \omega = -2 a_1\, e_1^* \wedge e_4^*
            - a_2\, e_2^* \wedge e_4^*
            - (a_2 + a_3)\, e_3^* \wedge e_4^*\]
%$d \omega = -2 a_1\, e_1^* \wedge e_4^*
%            - a_2\, e_2^* \wedge e_4^*
%            - (a_2 + a_3)\, e_3^* \wedge e_4^*$,
and~\eqref{Z2(g)andB2(g)} implies
\[
    B^2(\Cc, \K)
     = \s\{e_1^* \wedge e_4^*,\,
           e_2^* \wedge e_4^*,\,
           e_3^* \wedge e_4^*\}.
\]
Using~\eqref{second-exterior-derivative}, one checks that
$d(e_1^* \wedge e_4^*) = d(e_2^* \wedge e_4^*)
 = d(e_3^* \wedge e_4^*) = 0$, while the other three basis
2-forms have nonzero exterior
differentials. Hence
\[
    Z^2(\Cc, \K)
     = \s\{e_1^* \wedge e_4^*,\,
           e_2^* \wedge e_4^*,\,
           e_3^* \wedge e_4^*\}
     = B^2(\Cc, \K),
\]
so $\Cc$ is cocomplete.
\end{example}

We have used Algorithm~\ref{alg2} to determine all complex and real cocomplete Lie algebras of dimension at most 4, and the result is recorded in Table~\ref{tab-cocompleteLA} below. Throughout the table, the displayed parameter conditions ensure cocompleteness, the underlying Lie algebras themselves (for arbitrary parameter values) appear in~\cite{SW14}.

%-------Table 2-------------
\begin{table}[h]
    \centering
    \caption{Complex and real cocomplete Lie algebras of dimension $\leq 4$}\label{tab-cocompleteLA}
    \begin{tabular}{c p{7.8cm} c}
        \hline Lie algebras
            & \centering Non-zero Lie brackets, conditions of parameters
            & References \\
        \hline
            $\aff$ & $[e_1, e_2] = e_2$ & \\
        \hline
            $\mathfrak{sl}_2(\K)$
              & $[e_1, e_2] = e_2$, $[e_1, e_3] = -e_3$, $[e_2, e_3] = e_1$ & \\
            $\mathfrak{so}_3(\R)$
              & $[e_1, e_2] = e_3$, $[e_1, e_3] = -e_2$, $[e_2, e_3] = e_1$ & \\
            $\Sf_{3,1}$
              & $[e_1, e_3] = e_1$, $[e_2, e_3] = a e_2$ \; ($a \neq 0, -1$)
              & \cite[Sec.~16.4]{SW14} \\
            $\Sf_{3,2}$
              & $[e_1, e_3] = e_1$, $[e_2, e_3] = e_1 + e_2$
              & \cite[Sec.~16.4]{SW14} \\
            $\Sf_{3,3}\,(\text{over } \R)$
              & $[e_1, e_3] = a e_1 - e_2$, $[e_2, e_3] = e_1 + a e_2$ \; ($a \neq 0$)
              & \cite[Sec.~16.4]{SW14} \\
        \hline
            $\mathfrak{sl}_2(\K) \LO \K$
              & $[e_1, e_2] = e_2$, $[e_1, e_3] = -e_3$, $[e_2, e_3] = e_1$ & \\
            $\mathfrak{so}_3(\R) \LO \R$
              & $[e_1, e_2] = e_3$, $[e_1, e_3] = -e_2$, $[e_2, e_3] = e_1$ & \\
            $\Sf_{4,2}$
              & $[e_1, e_4] = e_1$, $[e_2, e_4] = e_1 + e_2$,
              
              $[e_3, e_4] = e_2 + e_3$
              & \cite[Sec.~17.2]{SW14} \\
            $\Sf_{4,3}$
              & $[e_1, e_4] = e_1$, $[e_2, e_4] = a e_2$, $[e_3, e_4] = b e_3$ \; ($a, b \neq 0, -1$; $a + b \neq 0$)
              & \cite[Sec.~17.2]{SW14} \\
            $\Sf_{4,4}$
              & $[e_1, e_4] = e_1$, $[e_2, e_4] = e_1 + e_2$, $[e_3, e_4] = a e_3$ \; ($a \neq 0, -1$)
              & \cite[Sec.~17.2]{SW14} \\
            $\Sf_{4,5}\,(\text{over } \R)$
              & $[e_1, e_4] = a e_1$, $[e_2, e_4] = b e_2 - e_3$, 
              
              $[e_3, e_4] = e_2 + b e_3$ \; ($a, b \neq 0$)
              & \cite[Sec.~17.2]{SW14} \\
            $\Sf_{4,8}$
              & $[e_1, e_4] = (1 + a) e_1$, $[e_2, e_4] = e_2$, 
              
              $[e_3, e_4] = a e_3$ \; ($a \notin \{0, -1, -\tfrac{1}{2}, -2\}$)
              & \cite[Sec.~17.3]{SW14} \\
            $\Sf_{4,9}\,(\text{over } \R)$
              & $[e_1, e_4] = 2 a e_1$, $[e_2, e_4] = a e_2 - e_3$, 
              
              $[e_3, e_4] = e_2 + a e_3$ \; ($a \neq 0$)
              & \cite[Sec.~17.3]{SW14} \\
            $\Sf_{4,10}$
              & $[e_1, e_4] = 2 e_1$, $[e_2, e_4] = e_2$, $[e_3, e_4] = e_2 + e_3$
              & \cite[Sec.~17.3]{SW14} \\
        \hline
    \end{tabular}
\end{table}

%----Subsection 4.4: Almost abelian cocomplete----
\subsection{Almost abelian cocomplete Lie algebras: Proof of Theorem~\ref{thm3}}\label{subsec:4-4}

Based on the criterion $H^2(\Cc, \K) = 0$ and the structural foundation from the previous subsections, we now focus on the class of almost abelian Lie algebras. The choice of this class is motivated not only by the fact that it represents the most natural extension of purely abelian structures, but also because the
existence of an abelian ideal of codimension one allows us to translate the abstract cohomological problem into an explicit linear algebra problem. We first recall this class.

%--Definition 4.12: Almost abelian Lie algebra (AALA)--
\begin{definition}[{Almost abelian Lie algebras~\cite[Definition~1]{Ave22}}]\label{D341}
A non-abelian Lie algebra $\Cc$ is called \emph{almost abelian} (or, for short, an \emph{AALA}) if it contains an abelian subalgebra of codimension 1. When $\Cc$ is an $n$-dimensional AALA ($1 < n \in \N$), it is also called an \emph{$n$-AALA}.
\end{definition}

%----Remark 4.13---
\begin{remark}\label{R342} We have the following remarks.
\begin{enumerate} 
    \item By~\cite[Proposition~3.1]{BC12}, the abelian subalgebra     in Definition~\ref{D341}, if any, is in fact an ideal of $\Cc$. Moreover, every AALA $\Cc$ is 2-step solvable, i.e.,
    $[\Cc, \Cc]$ is abelian.

    \item When $\Cc$ is an $(n+1)$-AALA, we can always view     $\K^n$ as a 1-codimensional ideal of $\Cc$. Therefore, $\Cc$ is a 1-dimensional nontrivial extension of $\K^n$ by some
    $D \in \Der(\K^n) \equiv \mathfrak{gl}_n(\K)$, i.e.,
    $\Cc = \K^n \oplus_D \K e_0$ with
    $e_0 \in \Cc \setminus \K^n$. Equivalently,
    $\Cc \equiv \K^n \oplus \K e_0$ with bracket
    $[e_0, v] \co D(v) \in \K^n$ for every $v \in \K^n$.
\end{enumerate}
\end{remark}

As mentioned at the end of Section~\ref{sec1}, we now prove
Theorem~\ref{thm3}, which gives a necessary and sufficient condition for an $(n+1)$-AALA $\Cc = \K^n \oplus_D \K e_0$ to be cocomplete. We begin with the following lemma.

%---Lemma 4.14---
\begin{lemma}\label{L343}
If $\Cc = \K^n \oplus_D \K e_0$ is an AALA, then
\[
    H^2(\Cc, \K) \cong
       \Lambda^2_D((\K^n)^*)
       \oplus \Lambda^1((\K^n)^*) / \im D^T.
\]
Here $D^T \in \Der((\K^n)^*)$ is the transpose of $D$, and
$\Lambda^2_D((\K^n)^*)$ denotes the subspace of 2-forms
$\omega \in \Lambda^2((\K^n)^*)$ such that
$D^T \cdot \omega \equiv 0$, where
\[
    (D^T \cdot \omega)(v_1, v_2)
     \co \omega(D(v_1), v_2) + \omega(v_1, D(v_2)),
    \quad v_1, v_2 \in \K^n.
\]
\end{lemma}

\begin{proof}
Since $\K^n$ has the canonical basis $(e_1, \ldots, e_n)$, $\Cc = \K^n \oplus \K e_0 \equiv \K^{n+1}$ admits the basis $B \co (e_0, e_1, \ldots, e_n)$. Let $(e_1^*, \ldots, e_n^*)$ be the dual basis of $(e_1, \ldots, e_n)$. Then $\Cc^* \cong (\K^{n+1})^*$ admits the dual basis $B^* \co (e_0^*, e_1^*, \ldots, e_n^*)$. For any 1-form $\eta \in \Lambda^1((\K^n)^*)$, the Hochschild--Serre formula~\cite[p.~592]{HS53} gives the 2-coboundary of $\eta$ as
\[
    d \eta = -e_0^* \wedge D^T \eta
          \in B^2(\Cc, \K) \subset C^2(\Cc, \K).
\]
Hence the space of 2-coboundaries on $\Cc$ is
\[
    B^2(\Cc, \K) = \s\{e_0^*\} \wedge \im D^T.
\]
Now consider the space $C^2(\Cc, \K)$ of 2-cochains. Every
$\omega \in C^2(\Cc, \K)$ can be represented as
\[
    \omega = \xi + \alpha\, e_0^* \wedge \eta
          \in \Lambda^2((\K^n)^*)
              \oplus \s\{e_0^*\} \wedge \Lambda^1((\K^n)^*),
    \quad \alpha \in \K.
\]
For $\xi \in \Lambda^2((\K^n)^*)$, the Hochschild--Serre formula yields
\begin{align*}
    d \xi(e_0, e_i, e_j)
     &= -\bigl(\xi(D(e_i), e_j) + \xi([e_i, e_j], e_0)
                + \xi(-D(e_j), e_i)\bigr) \\
     &= -\bigl(\xi(D(e_i), e_j) + \xi(e_i, D(e_j))\bigr) \\
     &= -(D^T \cdot \xi)(e_i, e_j),
        \quad 1 \leq i < j \leq n,
\end{align*}
i.e., $d \xi(e_0, \cdot, \cdot)
       = -(D^T \cdot \xi)(\cdot, \cdot)$.
Since $d \omega = d \xi + \alpha\, d(e_0^* \wedge \eta) = d \xi$,
\[
    d \omega = 0 \iff d \xi = 0
              \iff D^T \cdot \xi = 0
              \iff \xi \in \Lambda^2_D((\K^n)^*).
\]
Hence the space of 2-cocycles on $\Cc$ is
\[
    Z^2(\Cc, \K)
     = \Lambda^2_D((\K^n)^*)
       \oplus \s\{e_0^*\} \wedge \Lambda^1((\K^n)^*).
\]
By definition,% of the second cohomology,
\[
    H^2(\Cc, \K)
     \co Z^2(\Cc, \K) / B^2(\Cc, \K)
     = \Lambda^2_D((\K^n)^*)
       \oplus \Lambda^1((\K^n)^*) / \im D^T. \qedhere
\]
\end{proof}

Now we will prove Theorem~\ref{thm3}, the final main result of the paper.

\begin{proof}[{\bf Proof of Theorem~\ref{thm3}}]
By Proposition~\ref{P321} and Lemma~\ref{L343}, for any AALA $\Cc = \K^n \oplus_D \K e_0$, it is cocomplete if and only if
\[
    \Lambda^2_D((\K^n)^*)
        \oplus \Lambda^1((\K^n)^*) / \im D^T = 0
    \iff \Lambda^2_D((\K^n)^*)
         = \Lambda^1((\K^n)^*) / \im D^T = 0.
\]
We will prove that the two summands above vanish if and only if conditions~\eqref{thm3-1} and~\eqref{thm3-2} of Theorem~\ref{thm3} hold.
\begin{itemize}
    \item For the second summand $\Lambda^1((\K^n)^*) / \im D^T$, we have
    \[
        \Lambda^1((\K^n)^*) / \im D^T = 0
         \iff \im D^T = (\K^n)^*
         \iff \im D = \K^n,
    \]
    i.e., $D \colon \K^n \to \K^n$ is a linear isomorphism. Hence
    \begin{equation}\label{4-1}
        \Lambda^1((\K^n)^*) / \im D^T = 0
         \iff \bigl(\text{Item~\eqref{thm3-1}
                          of Theorem~\ref{thm3} holds}\bigr).
    \end{equation}

    \item For the first summand $\Lambda^2_D((\K^n)^*)$, suppose first that condition~\eqref{thm3-1} holds, i.e., $D \in \Au(\K^n) \equiv \GL_n(\K)$. Since $\bar{\K}$ is algebraically closed, $D_{\bar{\K}} \co D \otimes_\K \id_{\bar{\K}}$ is triangularizable, i.e., there exists a basis $(v_1, \ldots, v_n)$ of $\bar{\K}^n$ in which $D_{\bar{\K}}$ is upper-triangular with diagonal entries $\lambda_1, \ldots, \lambda_n$ (the eigenvalues of $D_{\bar{\K}}$ counted with multiplicity). In the dual basis $(v_1^*, \ldots, v_n^*)$, the family $(v_p^* \wedge v_q^* \mid 1 \leq p < q \leq n)$ is a basis of $\Lambda^2((\bar{\K}^n)^*)$, and a direct computation gives
    \[
        \big(D_{\bar{\K}}^T \cdot (v_p^* \wedge v_q^*)\big)(v_p, v_q)
         = \lambda_p + \lambda_q,
    \]
    while $\big(D_{\bar{\K}}^T \cdot (v_p^* \wedge v_q^*)\big)(v_i, v_j) = 0$ whenever $(i, j)$ exceeds $(p, q)$ in lexicographic order. Hence $D_{\bar{\K}}^T \cdot$ is upper-triangular in this basis with diagonal entries $\lambda_p + \lambda_q$ for $p < q$; in particular, its eigenvalues are exactly $\{\lambda_p + \lambda_q \mid 1 \leq p < q \leq n\}$.

    It follows that $D_{\bar{\K}}^T \cdot$ is injective on $\Lambda^2((\bar{\K}^n)^*)$. Equivalently, $\Lambda^2_{D_{\bar{\K}}}((\bar{\K}^n)^*) = 0$ if and only if $\lambda_p + \lambda_q \neq 0$ for all $1 \leq p < q \leq n$. Since
    \[
        \Lambda^2_{D_{\bar{\K}}}((\bar{\K}^n)^*)
         \cong \Lambda^2_D((\K^n)^*) \otimes_\K \bar{\K},
    \]
    we conclude
    \begin{equation}\label{4-2}
        \Lambda^2_D((\K^n)^*) = 0
         \iff \bigl(\text{Item~\eqref{thm3-2}
                          of Theorem~\ref{thm3} holds}\bigr).
    \end{equation}
\end{itemize}
Combining~\eqref{4-1} with~\eqref{4-2} yields the conclusion of Theorem~\ref{thm3}.
\end{proof}

%----Subsection 4.5: Classification of AALA--- 
\subsection{Classification of almost abelian cocomplete Lie algebras}\label{subsec:4-5}

We now address the classification of almost abelian cocomplete Lie algebras. The classification of arbitrary AALAs $\Cc = \K^n \oplus_D \K e_0$ is well-known, being reduced to the classification of $D \in \mathfrak{gl}_n(\K)$ as follows.

%---Proposition 4.15-----
\begin{proposition}[{Avetisyan~\cite[Proposition~10]{Ave22}}]\label{P351}
Two AALAs 
\[\Cc = \K^n \oplus_D \K e_0 \quad \text{and} \quad \Cc' = \K^n \oplus_{D'} \K e_0'\]
%$\Cc = \K^n \oplus_D \K e_0$ and
%$\Cc' = \K^n \oplus_{D'} \K e_0'$ 
are isomorphic if and only if %there exist 
$\exists \, \alpha \in \K \setminus \{0\}$, % and
$\exists \, \Phi \in \Au(\K^n) \equiv \GL_n(\K)$ such that
$\alpha D' = \Phi^{-1} D \Phi$.
\end{proposition}

The operators $D$ and $D'$ in Proposition~\ref{P351} are then called \emph{proportionally similar}. We recall the precise definition below.

%---Definition 4.16----
\begin{definition}[{Proportional similarity~\cite[Definition~2.1]{Le}}]\label{D352}
Two operators $D, D' \in \mathfrak{gl}(\K^n)$ (resp., two matrices $M, N \in \M_n(\K)$) are \emph{proportionally similar}, denoted by $D \sim_p D'$ (resp., $M \sim_p N$), if there exist $\alpha \in \K \setminus \{0\}$ and $\Phi \in \Au(\K^n) \equiv \GL_n(\K)$ (resp., $P \in \GL_n(\K)$) such that $\alpha D' = \Phi^{-1} D \Phi$ (resp., $\alpha N = P^{-1} M P$). We then also say $D \sim_p D'$ (resp., $M \sim_p N$) via $\alpha$ and $\Phi$ (resp., $P$).
\end{definition}

%---Remark 4.17---
\begin{remark}\label{R353}
The relation $\sim_p$ is an equivalence relation on $\M_n(\K)$. Moreover, if $M \sim_p N$ via $\alpha \in \K \setminus \{0\}$, then the Jordan canonical forms of $\alpha M$ and $N$ coincide (sizes of the Jordan blocks are equal, while eigenvalues are scaled by $\alpha$). Therefore, the quotient set $\M_n(\K)/{\sim_p}$ can be determined using Jordan canonical forms.
\end{remark}

By Remark~\ref{R353}, Proposition~\ref{P351} reduces the
classification of $(n+1)$-dimensional almost abelian Lie algebras to the determination of $\M_n(\K)/{\sim_p}$. Combined with Theorem~\ref{thm3}, this gives the following.

%---Proposition 4.18---
\begin{proposition}\label{classification-aac}
Let $\Cc = \K^n \oplus_D \K e_0$ and
$\Cc' = \K^n \oplus_{D'} \K e_0$ be almost abelian cocomplete Lie algebras, where $D, D' \in \Au(\K^n) \equiv \GL_n(\K)$ satisfy condition~\eqref{thm3-2} of Theorem~\ref{thm3}. Then $\Cc \cong \Cc'$ if and only if $D \sim_p D'$.
\end{proposition}

In particular, if $\GL_n^0(\K) \subset \GL_n(\K)$ denotes the subset of invertible matrices satisfying
condition~\eqref{thm3-2} of Theorem~\ref{thm3} (i.e., $D$ has no pair of eigenvalues summing to zero), then
Proposition~\ref{classification-aac} yields the following.

%----Proposition 4.19---
\begin{proposition}\label{P354}
The classification of $(n+1)$-dimensional almost abelian
cocomplete Lie algebras $\Cc = \K^n \oplus_D \K e_0$ is
equivalent to the classification of $D \in \GL_n^0(\K)$ up to proportional similarity. Equivalently, one has to determine the quotient set $\GL_n^0(\K)/{\sim_p}$.
\end{proposition}

Proposition~\ref{P354} provides the theoretical foundation for Algorithm~\ref{alg3}, which classifies almost abelian cocomplete Lie algebras. To obtain the classification up to isomorphism, we use the additional algorithms in~\cite{NLV25}.

%---Algorithm 3---
\begin{algorithm}[h]
    \KwIn{A positive integer $n$}
    \KwOut{A list $\Sc$ of non-isomorphic $(n+1)$-dimensional almost abelian cocomplete Lie algebras}
    $\Sc \co \{\,\}$\;
    $\Cc_D \co \K^n \oplus_D \K e_0$ with $D \in \GL_n(\K)$\;
    $\GL_n^0(\K) \co
       \{D \in \GL_n(\K) \text{ satisfies
         condition~\eqref{thm3-2} of Theorem~\ref{thm3}}\}$\;
    $\mathcal{D} \co \GL_n^0(\K)/{\sim_p}$\;
    \For{$D \in \mathcal{D}$}{
        Isomorphism verification for $\Cc_D$ by algorithms
        in~\cite{NLV25}\;
        $\mathcal{I} \co
          \{\text{Isomorphism classes of } \Cc_D\}$\;
        $\Sc \co \Sc \cup \mathcal{I}$\;
    }
    Return $\Sc$.
    \caption{Classifying almost abelian cocomplete Lie algebras}\label{alg3}
\end{algorithm}

%---Example 4.20---
\begin{example}\label{Ex-classif-3dim}
Let $\Cc = \K^2 \oplus_D \K e_0$ be a 3-dimensional almost abelian cocomplete Lie algebra with
$D \in \Au(\K^2) \equiv \GL_2(\K)$. By
condition~\eqref{thm3-2} of Theorem~\ref{thm3}, $D$ has no pair of nonzero eigenvalues summing to zero. Over $\K = \C$, the possible Jordan canonical forms of $D$ are
\[
    J_1 \co \begin{bmatrix} \lambda_1 & 0 \\ 0 & \lambda_2 \end{bmatrix}
        \;(\lambda_1 \lambda_2 \neq 0, \;
           \lambda_1 \neq -\lambda_2)
    \quad \text{or} \quad
    J_2 \co \begin{bmatrix} \lambda & 1 \\ 0 & \lambda \end{bmatrix}
        \;(\lambda \neq 0).
\]
By scaling,
\[
    J_1 \sim_p \begin{bmatrix} 1 & 0 \\ 0 & \lambda \end{bmatrix}
        \co D_1
        \;\Bigl(\lambda = \tfrac{\lambda_2}{\lambda_1} \neq 0, -1\Bigr)
    \quad \text{and} \quad
    J_2 \sim_p \begin{bmatrix} 1 & 1 \\ 0 & 1 \end{bmatrix}
        \co D_2.
\]
When $\K = \R$, the real Jordan canonical form
\[
    J_3 \co \begin{bmatrix} \lambda & 1 \\ -1 & \lambda \end{bmatrix}
        \; (\lambda \neq 0)
\]
also occurs. Hence, 
\[\GL_2^0(\C)/{\sim_p} = \{D_1, D_2\}, \quad \GL_2^0(\R)/{\sim_p} = \{D_1, D_2, J_3\}.\]
%$\GL_2^0(\C)/{\sim_p} = \{D_1, D_2\}$, while
%$\GL_2^0(\R)/{\sim_p} = \{D_1, D_2, J_3\}$. 
These quotient sets give the following 3-dimensional almost abelian cocomplete Lie algebras in the basis $(e_0, e_1, e_2)$ as follows:
\[
    \begin{array}{l l l}
        \Cc_{3.1}^\lambda \colon
         & [e_0, e_1] = e_1,\; [e_0, e_2] = \lambda e_2
         & (\lambda \neq 0, -1),\\
        \Cc_{3.2} \colon
         & [e_0, e_1] = e_1,\; [e_0, e_2] = e_1 + e_2,\\
        \rc_{3.3}^\lambda \colon
         & [e_0, e_1] = \lambda e_1 - e_2,\;
           [e_0, e_2] = e_1 + \lambda e_2
         & (\lambda \neq 0).
    \end{array}
\]
Here, $\rc$ indicates a Lie algebra that exists only over $\R$, the others are valid over both $\R$ and $\C$. Using algorithms in~\cite{NLV25}, one finds that
$\Cc_{3.1}^\lambda \cong \Cc_{3.1}^{1/\lambda}$ and
$\rc_{3.3}^\lambda \cong \rc_{3.3}^{-\lambda}$ via the
isomorphisms
\[
    \begin{bmatrix}
        \lambda & 0 & 0 \\
        0 & 0 & 1 \\
        0 & 1 & 0
    \end{bmatrix}
    \quad \text{and} \quad
    \begin{bmatrix}
        -1 & 0 & 0 \\
        0 & 0 & 1 \\
        0 & 1 & 0
    \end{bmatrix},
\]
respectively. %; the notations $\lambda \equiv \tfrac{1}{\lambda}$ and $\lambda \equiv -\lambda$ are used to indicate these identifications. 
In summary, we have classified 3-dimensional complex and real almost abelian cocomplete Lie algebras, up to isomorphism.
\end{example}

We have used Algorithm~\ref{alg3} to classify, up to isomorphism, the complex and real almost abelian cocomplete Lie algebras of dimension at most 4. The results are presented in Table~\ref{tab-aaCLA}, in which ``$\equiv$'' indicates that the
corresponding parameters give rise to isomorphic Lie algebras, and the absence of ``$\equiv$'' means that the parameter is optimal in the sense that distinct parameters yield non-isomorphic Lie algebras.

%--------------Table 3 ----
\begin{table}[!h]
    \centering
    \caption{Complex and real almost abelian cocomplete Lie algebras of dimension $\leq 4$}\label{tab-aaCLA}
    \begin{tabular}{c p{10.5cm}}
        \hline Lie algebras
            & \centering Non-zero Lie brackets, conditions of parameters
              and isomorphisms \tabularnewline
        \hline
            $\aff$ & $[e_0, e_1] = e_1$ \\
        \hline
            $\Cc_{3.1}^\lambda$
              & $[e_0, e_1] = e_1$, $[e_0, e_2] = \lambda e_2$
                \;\big($\lambda \neq 0, -1$; $\lambda \equiv \tfrac{1}{\lambda}$\big) \\
            $\Cc_{3.2}$
              & $[e_0, e_1] = e_1$, $[e_0, e_2] = e_1 + e_2$ \\
            $\rc_{3.3}^\lambda$
              & $[e_0, e_1] = \lambda e_1 - e_2$,
                $[e_0, e_2] = e_1 + \lambda e_2$
                \;($\lambda \neq 0$; $\lambda \equiv -\lambda$) \\
        \hline
            \multirow{2.5}{*}{$\Cc_{4.1}^{\lambda\mu}$}
              & $[e_0, e_1] = e_1$, $[e_0, e_2] = \lambda e_2$,
                $[e_0, e_3] = \mu e_3$
                \;
                
                ($\lambda, \mu \neq 0, -1$;
                   $\lambda + \mu \neq 0$) \\
              & $\bigl((\lambda, \mu)
                  \equiv (\mu, \lambda)
                  \equiv (\tfrac{\lambda}{\mu}, \tfrac{1}{\mu})
                  \equiv (\tfrac{\mu}{\lambda}, \tfrac{1}{\lambda})
                  \equiv (\tfrac{1}{\mu}, \tfrac{\lambda}{\mu})
                  \equiv (\tfrac{1}{\lambda}, \tfrac{\mu}{\lambda})
                \bigr)$ \\
            $\Cc_{4.2}^\lambda$
              & $[e_0, e_1] = e_1$, $[e_0, e_2] = \lambda e_2$,
                $[e_0, e_3] = e_2 + \lambda e_3$
                \;($\lambda \neq 0, -1$) \\
            $\rc_{4.3}^{\lambda\mu}$
              & $[e_0, e_1] = \lambda e_1$, $[e_0, e_2] = \mu e_2 - e_3$,
                $[e_0, e_3] = e_2 + \mu e_3$
                \;
                
                ($\lambda, \mu \neq 0$) \\
            $\Cc_{4.4}$
              & $[e_0, e_1] = e_1$, $[e_0, e_2] = e_1 + e_2$,
                $[e_0, e_3] = e_2 + e_3$ \\
            \hline
    \end{tabular}
\end{table}

%--------- Concluding Remark-------------
\section*{Concluding remarks}
\addcontentsline{toc}{section}{Concluding remarks}

In this paper, we have investigated the ``dual manifestations'' of injective and projective properties within the category of Lie algebras, leading to the formal introduction of cocomplete Lie algebras. Our primary contributions are summarized as follows:

\begin{itemize}
    \item First, we established that complete Lie algebras play an injective-type role in the category of Lie algebras: a Lie algebra trivially splits every extension by it if and only if it is complete (Theorem~\ref{thm1}). The four equivalent conditions in Theorem~\ref{thm1} clarify this role from complementary internal, lifting, and external-splitting viewpoints. Nevertheless, as Example~\ref{Ex-counter-cat-inj} shows, this injective-type property is strictly weaker than categorical injectivity.
    
    \item Second, we showed that no global projective analogue exists over the entire category of Lie algebras (Theorem~\ref{thm2}). This obstruction motivated the restriction to the natural and extensively studied class of central extensions.
    
    \item Third, to restore the symmetry of the ``injective-projective duality'' in this restricted setting, we introduced the notion of \emph{cocomplete} Lie algebras (Definition~\ref{D311}). Theorem~\ref{P-lifting-cocomplete} characterizes cocompleteness through four equivalent conditions: the trivial splitting of every central extension, a direct-summand property for central kernels, and a homomorphism-lifting property dual to that of complete Lie algebras. The latter mirrors, within the central setting, the lifting property of projective objects in the category of modules. Moreover, cocompleteness admits an explicit cohomological characterization: a Lie algebra $\Cc$ is cocomplete if and only if its second cohomology with trivial coefficients vanishes, $H^2(\Cc, \K) = 0$ (Proposition~\ref{P321}); in particular, semisimple Lie algebras are cocomplete.
  
    \item Finally, Theorem~\ref{thm3} provides an explicit spectral criterion for an almost abelian Lie algebra to be cocomplete, and Proposition~\ref{P354} reduces their classification %of these Lie algebras 
    to that of invertible matrices up to proportional similarity. Three corresponding algorithms support these characterizations and yield tabulations of the complete, the cocomplete, and the almost abelian cocomplete Lie algebras of dimension at most 4.   
\end{itemize}

Our future research will focus on the classification of complete and cocomplete Lie algebras in low dimensions or under further structural conditions. In particular, we intend to explicitly describe finite-dimensional cocomplete Lie algebras with maximal commutative ideals of codimension 2.

\section*{Acknowledgements}

This research is funded by Vietnam National University Ho Chi Minh City (VNU-HCM) under grant number B-2026-34-06.

%--------- References-----------------------


\addcontentsline{toc}{section}{References}
\begin{thebibliography}{99}

\bibitem{Ave22}Z. Avetisyan, The structure of almost Abelian Lie algebras, Internat. J. Math. 33 (8) (2022) 2250057 (26 pages).

\bibitem{BC12}D. Burde, M. Ceballos, Abelian ideals of maximal dimension for solvable Lie algebras, J. Lie Theory 22 (3) (2012) 741--756.

%\bibitem{CCL99}S. S. Chern, W. H. Chen, K. S. Lam, Lectures on Differential Geometry, Series on University Mathematics Vol. 1, World Scientific, Singapore, 1999.

%\bibitem{Che44}C. Chevalley, On groups of automorphism of Lie groups,	Proc. Natl. Acad. Sci. USA 30 (1944) 274--275.

\bibitem{Che-Eil}C. Chevalley, S. Eilenberg, Cohomology theory of Lie groups and Lie algebras, Trans. Amer. Math. Soc. 63 (1) (1948) 85--124.

\bibitem{Dix55}J. Dixmier, Cohomologie des alg\`{e}bres de Lie nilpotentes, Acta Sci. Math. (Szeged) 16 (1955) 246--250.

%\bibitem{Gei76}A. G. Gein, Semimodular Lie algebras, Sib. Math. J. 17 (1976) 189--193.
	
%\bibitem{Gra05}W. A. de Graaf, Classification of solvable Lie algebras, Exp. Math. 14 (1) (2005) 15--25.

\bibitem{HS53} G. Hochschild, J-P. Serre, Cohomology of Lie algebras, Ann. of Math. 57 (3) (1953) 591--603.

\bibitem{Jac62}N. Jacobson, Lie Algebras, Wiley, New York, 1962.

\bibitem{Kas82}F. Kasch, Modules and Rings, London Mathematical Society Monograph No. 17, Academic Press, London -- New York, 1982.

\bibitem{KBK97}E. A. de Kerf, G. G. A. B\"auerle, A. P. E. ten Kroode, Lie Algebras, part 2: Finite and Infinite Dimensional Lie Algebras and Applications in Physics, Studies in Mathematical Physics Vol. 7, North-Holland, The Netherlands, 1997.

%\bibitem{Kol65}B. Kolman, Semi-modular Lie algebras, J. Sci. Hiroshima Univ., Ser. A-1 29 (2) (1965) 149--163.

\bibitem{Le}V. A. Le, H. T. T. Cao, H. Q. Duong, T. A. Nguyen, T. N. Vo, On the problem of classifying solvable Lie algebras having small codimensional derived algebras, Comm. Algebra 50 (9) (2022) 3775--3793.

\bibitem{NLV25}T. A. Nguyen, V. A. Le, T. N. Vo, Testing isomorphism of complex Lie algebras, Comm. Algebra 54 (4) (2026) 1320--1338.%, doi.org/10.1080/00927872.2025.2553026

\bibitem{Sch08}M. Schottenloher, A Mathematical Introduction to Conformal Field Theory, Lecture Notes in Physics Vol. 759, Springer-Verlag, Berlin -- Heidelberg, 2008.

\bibitem{Serre}J. P. Serre, Lie Algebras and Lie Groups (1964 Lectures given at Harvard University), Springer-Verlag, Berlin -- Heidelberg, 1992.

\bibitem{SW14}L. \v Snobl, P. Winternitz, Classification and Identification of Lie Algebras, CRM Monograph Series Vol. 33, American Mathematical Society, Providence, 2014.

%\bibitem{To67}L. \v S. T\^og\^o, Outer derivations of Lie algebras, Trans. Amer. Math. Soc. 128 (2) (1967) 264--276.

\bibitem{Weibel}C. A. Weibel,  An Introduction to Homological Algebra, Cambridge Studies in Advanced Mathematics 38, Cambridge University Press, UK, 1994. 



\end{thebibliography}
\end{document}